\documentclass[aap,MSNbibl,seceqn,citesort,dvips]{arximspdf}
\doi{10.1214/11-AAP797}
\volume{22}
\issue{4}
\pubyear{2012}
\firstpage{1411}
\lastpage{1449}

\makeatletter
\newcommand{\iint}{\int\!\!\int}

\newcommand{\veps}{\varepsilon}
\newcommand{\rmd}{\mathrm{d}}
\newcommand{\rmi}{\mathrm{i}}
\newcommand{\wh}{\widehat}
\newcommand{\whH}{\widehat{H}}
\newcommand{\whV}{\widehat{V}}
\newcommand{\whL}{\widehat{L}}
\newcommand{\whX}{\widehat{X}}
\newcommand{\whk}{\widehat{\kappa}}
\newcommand{\pibar}{\overline{\Pi}}
\newcommand{\Vbar}{\overline{V}}
\newcommand{\Xbar}{\overline{X}}
\newcommand{\Ybar}{\overline{Y}}
\newcommand{\Wbar}{\overline{W}}
\newcommand{\Zbar}{\overline{Z}}
\newcommand{\R}{\mathbb{R}}
\newcommand{\tux}{\tau(u-x)}
\newcommand{\tu}{\tau(u)}
\newcommand{\too}{\tau(0)}
\newcommand{\tD}{\tau(\Delta)}
\newcommand{\tW}{\tau^W}
\newcommand{\XD}{X\circ\theta_{\tau(u-x)}}
\newcommand{\X}{X_{[0,\tau(u-x))}}
\newcommand{\Xt}{X_{[0,t)}}
\newcommand{\wt}{w_{[0,t)}}
\newcommand{\wta}{w_{[0,\rho)}}
\newcommand{\ca}{\beta_1}
\newcommand{\cb}{\beta_2}
\newcommand{\cab}{\beta}
\newcommand{\Pu}{P^{(u)}}
\newcommand{\Eu}{E^{(u)}}
\newcommand{\tZ}{\tilde Z}
\newcommand{\ga}{\alpha}
\newcommand{\gl}{\lambda}
\newcommand{\gk}{\kappa}
\newcommand{\gz}{\zeta}
\newcommand{\gt}{\theta}
\newcommand{\tmu}{\tilde\mu}
\newcommand{\tnu}{\tilde\nu}

\newcommand{\eqdr}{\stackrel{\mathrm{D}}{=}}

\newtheorem{theorem}{Theorem}[section]
\newtheorem{lemma}{Lemma}[section]
\newtheorem{prop}{Proposition}[section]
\newproclaim{remark}{Remark}[section]

\makeatother

\begin{document}
\begin{frontmatter}

\title{Path decomposition of ruinous behavior
for a~general L\'{e}vy insurance risk process}
\runtitle{Path decomposition of a L\'{e}vy insurance risk process}

\begin{aug}
\author[A]{\fnms{Philip S.} \snm{Griffin}\corref{}\ead[label=e1]{psgriffi@syr.edu}}
\and
\author[B]{\fnms{Ross A.} \snm{Maller}\thanksref{t1}\ead[label=e2]{Ross.Maller@anu.edu.au}}
\runauthor{P. S. Griffin and R. A. Maller}
\affiliation{Syracuse University and Australian National University}
\address[A]{Department of Mathematics\\
Syracuse University\\
Syracuse, New York 13244-1150\\
USA\\
\printead{e1}}
\address[B]{Centre for Financial Mathematics\\
MSI and School of Finance\\
\quad and Applied Statistics\\
Australian National University\\
Canberra, ACT 0200\\
Australia\\
\printead{e2}}
\end{aug}

\thankstext{t1}{Supported in part by ARC Grant DP1092502.}

\received{\smonth{5} \syear{2010}}
\revised{\smonth{6} \syear{2011}}

%
\begin{abstract}
We analyze the general L\'evy insurance risk process for L\'evy
measures in the convolution equivalence class $\mathcal{S}^{(\alpha)}$,
$\alpha> 0$, via a new kind of path decomposition. This yields a very
general functional limit theorem as the initial reserve level $u\to
\infty$, and a host of new results for functionals of interest in
insurance risk. Particular emphasis is placed on the time to ruin,
which is shown to have a proper limiting distribution, as $u\to
\infty$, conditional on ruin occurring under our assumptions. Existing
asymptotic results under the $\mathcal{S}^{(\alpha)}$ assumption are
synthesized and extended, and proofs are much simplified, by
comparison with previous methods specific to the convolution
equivalence analyses. Additionally, limiting expressions for penalty
functions of the type introduced into actuarial mathematics by Gerber
and Shiu are derived as straightforward applications of our main
results.
\end{abstract}

%

\pdfkeywords{60G51, 60F17, 91B30,
62P05, Levy insurance risk process,
convolution equivalence,
time to ruin, overshoot,
expected discounted penalty function}

\end{frontmatter}


\section{Introduction}\label{s1}

Let $X=\{X_{t}\dvtx t \geq0 \}$,
$X_0=0$, be a L\'{e}vy process
defined on $(\Omega, {\mathcal F}, P)$,
with triplet $(\gamma, \sigma^2, \Pi_X)$,
$\Pi_{X}$ being the L\'{e}vy measure of~$X$.
Thus the characteristic function of $X$
is given by the L\'{e}vy--Khintchine
representation, $Ee^{i\theta X_{t}} = e^{t \Psi_X(\theta)}$,
where
\[
\Psi_X(\theta) =
\rmi\theta\gamma- \sigma^2\theta^2/2+
\int_{\R\setminus\{0\}}\bigl(e^{\rmi\theta x}-1-
\rmi\theta x \mathbf{1}_{\{|x|<1\}}\bigr)\Pi_X(\rmd x)\qquad
\mbox{for } \theta\in\mathbb{R}.
\]

We will be concerned with the case where $X_t\to-\infty$ a.s.
We have in mind an
insurance risk model with premiums and other income producing a
downward drift in $X$, while claims are represented by positive
jumps. Thus the process
$X$, called the \textit{claim surplus process}, represents the excess in
claims over premium.
We think of an insurance company starting with
an initial positive reserve
$u$, and ruin occurring
if this level is exceeded by~$X$.
We will refer to this as the \textit{general L\'evy insurance risk model}.
It is a~generalization of the classical \textit{Cram\'{e}r--Lundberg model},
which arises when the claim surplus process is taken to be
%
%
\begin{equation}\label{CL}
X_t=\sum_1^{N_t} U_i - rt,
\end{equation}
where $N_t$ is a Poisson process, $U_i> 0$ form an independent i.i.d.
sequence and $r>0$. Here $r$ represents the rate of premium inflow and
$U_i$ the size of the $i$th claim.
The general model allows for income other than through premium inflow
and a more realistic claims structure; see Section 2.7.1 of Kyprianou
\cite{kypbook}.
The assumption $X_t\to-\infty$ a.s. is a reflection of premiums being
set to avoid almost certain ruin for finite $u$.

The primary focus of this paper is on when and how ruin occurs for large
reserve levels, that is, as $u\to\infty$. Introduce
%
%
\begin{equation}\label{GXbar}
\Xbar_t=\sup_{0\le s\le t} X_s,\qquad
G_t= \sup\{0\le s\le t\dvtx X_s = \Xbar_{s}\}
\end{equation}
and
%
%
\begin{equation}\label{tau}
\tau(u)= \inf\{t \geq0 \dvtx X_{t}>u \}.
\end{equation}
(In cases where possible confusion might arise, we
will indicate the dependence on the process under consideration
by a superscript, as in $G_t^X$.) These variables play a central role
in fluctuation theory for L\'evy processes, and give rise to the main
variables of interest in insurance risk:

\begin{itemize}
\item Ruin time:
$ \tau(u)$,

\item Shortfall at ruin (overshoot): $X_{\tau(u)}-u$,

\item Surplus immediately prior to ruin (undershoot):
$u-X_{\tau(u)-}$,

\item Minimum surplus prior to ruin:
$u-\overline X_{\tau(u)-}$,

\item Time of minimum surplus prior to ruin: $G_{\tau(u)-}$,

\item Time remaining to ruin from the time of minimum
surplus: $\tau(u)-G_{\tau(u)-}$.
\end{itemize}



Our main interest is in the behavior of the process when ruin occurs,
that is, when $\tau(u)<\infty$.
Crucial questions, for example, are ``how long does it take
for ruin to occur?'' and ``what do the paths look like leading up to
ruin?''
We pay particular attention to these issues. We will exclude the
trivial case that $X$ is the negative of a subordinator, so $P(\tau
(u)<\infty)>0$ for finite $u;$ cf. (\ref{PK}) below. On the other hand,
the assumption $X_t\to-\infty$ a.s. implies
$P(\tau(u)<\infty)\to0$ as the initial level $u\to\infty$.
Consequently, it is convenient to introduce, by elementary means, a new
probability measure $P^{(u)}$ given by
\[
P^{(u)}( \cdot)=P\bigl(\cdot|\tau(u)<\infty\bigr),
\]
and to state our results as limit theorems
conditional on $\tau(u)<\infty$, that is, under~$P^{(u)}$.\vadjust{\goodbreak}

Some further background is useful to place our results in context.
The original work on the Cram\'{e}r--Lundberg model was done under the
\textit{Cram\'{e}r--Lundberg condition}
%
%
\begin{equation}\label{Cr}
E e^{\alpha X_1}=1 \qquad\mbox{for some } \alpha> 0,
\end{equation}
which among other things implies $X_t\to-\infty$ a.s.
Embrechts, Kl\"{u}ppelberg and Mikosch \cite{EKM}
call (\ref{Cr}) the \textit{small claims condition}. The major results
in this area include a large deviation estimate for the probability
of ruin
%
%
\begin{equation}\label{ruinprob}
e^{\alpha u}P\bigl(\tau(u)<\infty\bigr)\to C,
\end{equation}
where $C$ is a constant which can be identified, and $C>0$ if
%
%
\begin{equation}\label{Cr+}
E X_1e^{\alpha X_1}<\infty.
\end{equation}
In addition, the asymptotic behavior under $\Pu$ of several of the
variables listed above is known (see, e.g., \cite{A1} or \cite{EKM}).
The ruin estimate (\ref{ruinprob}) was extended to general L\'evy
insurance risk processes satisfying (\ref{Cr}) by Bertoin and Doney
\cite{BD94}.

A second regime under which the Cram\'{e}r--Lundberg model has been
studied is the \textit{subexponential} or
\textit{large claims case} (see Asmussen and Kl\"{u}ppel\-berg~\cite{AK}).
In this scenario, the claim size distribution
is subexponential and, roughly speaking,
ruin occurs solely due to the realization of one extremely large claim.


The small and large claims models each have their various
strengths and weaknesses.
A third, intermediate, regime was introduced recently in the general
model by Kl\"{u}ppelberg, Kyprianou and Maller \cite{kkm}.
To motivate this model,
observe that in the small claims case (\ref{Cr}) holds, while in
the large claims (subexponential) case
%
%
\begin{equation}\label{SubE}
E e^{\alpha X_1}=\infty\qquad\mbox{for all } \alpha> 0.
\end{equation}
Thus, to obtain a new model we must either consider processes whose
distributions satisfy (\ref{SubE}) and which are not subexponential
or processes which satisfy $E e^{\alpha X_1}<\infty$ for some $\alpha
> 0$
but for which (\ref{Cr}) fails.
It is the latter alternative that we will focus on.
Since $X_t\to-\infty$ a.s., it is easy to see that such processes
must satisfy that, for some $\alpha> 0$,
%
%
\begin{equation}\label{beta}
E e^{\alpha X_1}<1 \quad\mbox{and}\quad E e^{(\ga+\veps) X_1}=\infty
\qquad\mbox{for all } \veps> 0.
\end{equation}
For example, those with distribution tails of the form
%
%
\begin{equation}\label{p}
P(X_1>x)\sim\frac{e^{-\alpha x}}{x^p} \qquad\mbox{for } p>1
\end{equation}
satisfy (\ref{beta}).
A natural class of distributions which include those of the
form~(\ref{p}) is the class of convolution equivalent distributions of
index $\alpha$, which we now briefly describe.
As in \cite{kkm},
we will restrict ourselves to the nonlattice case, with\vadjust{\goodbreak}
the understanding that the alternative can be handled by obvious
modifications.
A~distribution $F$
on $ [0, \infty)$ with tail $\overline{F}=1-F$ belongs to
the \textit{class} $\mathcal{S}^{(\alpha)}$, $\alpha> 0$,
if $\overline{F}(u)>0$ for all $u>0$,
%
%
\begin{eqnarray}\label{Salph}
\lim_{u \to\infty}
\frac{\overline{F}(u+x)}{\overline{F}(u)} =
e^{-\alpha x}\qquad\mbox{for } x \in(-\infty, \infty)
\end{eqnarray}
and
%
%
\begin{equation}\label{S2}
\lim_{u \to\infty}
\frac{\overline{F}{}^{2*}(u)}{\overline{F}(u)}
\qquad\mbox{exists and is finite},
\end{equation}
where $F^{2*}=F * F$.
Distributions in $\mathcal{S}^{(\alpha)}$ are called
\textit{convolution equivalent} with index $\alpha$.
When $F\in\mathcal{S}^{(\alpha)}$, the limit in (\ref{S2})
must be of the form $ 2\delta_{\alpha}^F$,
where $\delta_{\alpha}^F:= \int_{[0, \infty)} e^{\alpha
x}F(\rmd x)$ is finite.
Much is known about the properties of such distributions. In
particular, the class is closed under tail equivalence,
that is, if $F\in\mathcal{S}^{(\alpha)}$ and $G$ is a
distribution function for which
\[
\lim_{u\to\infty}\frac{\overline{G}(u)}{\overline{F}(u)}=c
\qquad\mbox{for some } c\in(0,\infty),
\]
then $G\in\mathcal{S}^{(\alpha)}$.

Although the exponential distribution with parameter $\alpha$ is not
in $\mathcal{S}^{(\alpha)}$, distributions in $\mathcal{S}^{(\alpha
)}$ are
``near to exponential;'' for example, distributions with tails
comparable to $x^{-p}e^{-\alpha x}$, where $p>1$,
are in $\mathcal{S}^{(\alpha)}$.
The inverse Gaussian distributions, with appropriate choices of
parameters, form an important class of distributions which are
convolution equivalent. These in turn are a special case of the
tempered stable distributions, which have been the subject of
considerable recent activity.
For further examples and more on convolution equivalence
see \cite{C,EG,kl,Pakes} and~\cite{Pakes2}.

We can take the tail of any L\'evy measure, assumed
nonzero on some interval $(x_0,\infty)$, $x_0>0$, to be the
tail of a distribution function on $[0,\infty)$, after
renormalization.
With this convention, we say then that the measure (or its tail)
is in $\mathcal{S}^{(
\alpha
)}$ if this is true for the distribution with the
corresponding (renormalized) tail. The convolution equivalent model
introduced in \cite{kkm} is then one in which
%
%
\begin{equation}\label{c1}
\overline{\Pi}^+_{X} \in\mathcal{S}^{(\alpha)}
\quad\mbox{and}\quad
Ee^{\alpha X_1}< 1\qquad \mbox{for some } \alpha>0,
\end{equation}
%
where ${\Pi}^+_{X}$ is the restriction of ${\Pi}_{X}$ to $(0,\infty)$,
and $\pibar^+(x)=\Pi_X((x,\infty))$, \mbox{$x>0$}.
The condition $Ee^{\alpha X_1}< 1 $ implies
$e^{\alpha X_t}$ is a nonnegative supermartingale, from which it
follows immediately that $X_t\to-\infty$ a.s.
(This is also true when \mbox{$ Ee^{\alpha X_1}=1$}.)


By way of comparison with the small claims model, consider a one
parameter family of Cram\'{e}r--Lundberg models (\ref{CL}), in which
the claim size distribution \mbox{$U\in\mathcal{S}^{(\alpha)}$}. Let
\[
X^{(r)}_t=\sum_1^{N_t} U_i - rt,\qquad r\ge0,
\]
and set
\[
r_L = \frac{\ln(Ee^{\alpha X^{(0)}_1})}{\alpha}.
\]
Then $Ee^{\alpha X^{(r)}_1}=1$ if $r=r_L$ and
$Ee^{\alpha X^{(r)}_1}<1$ if $r>r_L$.
Thus the convolution equivalent models correspond to larger premium
rates (faster drift of $X$ to $-\infty$, lower probability of ruin),
than under the small claims condition~(\ref{Cr}).
In general, for any convolution equivalent model, there is an
associated model in which~(\ref{Cr}) holds, obtained by adding an
appropriate positive drift, which corresponds to decreasing the
premium rate. However, this change in premium rate leads to quite
different behavior in the two models.

Conditional on ruin occurring, the qualitative behavior of the claims
surplus process is very different in the convolution equivalent model
as opposed to either the small or large claims models.
In these latter two cases, the time to ruin, $\tau(u)$,
is of order $u$ as $u\to\infty$. In the small claims case,
under mild assumptions, there is a constant $b>0$ such that
\[
\frac{\tu}{u}\to b^{-1} \qquad\mbox{in } \Pu\mbox{ probability }
\]
and
\[
\sup_{t\in[0,1]}\biggl|\frac{X(t\tu)}{\tu}-bt\biggr|\to0
\qquad\mbox{in } \Pu\mbox{ probability},
\]
indicating that ruin occurs owing to the build up of small claims
which tend to cause $X$ to behave as though it had positive drift (see
\cite{A1} or \cite{EKM}).
In the subexponential case, the ruin time is again of
order $u$ (in distribution).
However, in this case the process evolves quite normally, that is,
like a sample path for which ruin does not occur, until a very large
claim suddenly causes ruin. This claim is so large that the
shortfall $X_{\tau(u)}-u{\stackrel{{P^{(u)}}}{\longrightarrow
}}\infty$
(see \cite{AK} or \cite{EKM}).

An obvious shortcoming of the small claims model is that it
does not allow for disasters, that is large jumps, which are observed
in real insurance data. On the other hand, the subexponential model is
very extreme and
uninformative in the sense that paths leading to ruin look quite normal
until suddenly
a large claim occurs, which results in
ruin with an arbitrarily large shortfall.

By contrast, the convolution equivalent model allows for disasters
to occur,
but they are not so ruinous as to be disproportionate in size relative
to the reserve level.
We will show that, in this model, asymptotically,
ruin occurs in finite time (in distribution), and for ruin to occur,
the claims surplus process must take a large jump from a neighborhood
of the origin to a neighborhood of $u$. This jump may result
in ruin, but if not, the process $X-u$ subsequently behaves like
$X$ conditioned to hit $(0,\infty)$. In either case,
the shortfall\vadjust{\goodbreak} at ruin converges in distribution to a finite random variable
as $u\to\infty$.
These results will follow from a path decomposition and asymptotic
analysis of the distribution of $X$,
conditional on ruin, in a way described below.
The idea of studying
ruin through a description of
the entire path leading up to ruin,
seems to have first appeared in Asmussen \cite{A1}, where the small
claims case for random walk is investigated. For work in the
subexponential case, see Asmussen and Kl\"{u}ppelberg \cite{AK}.




\section{Skorohod space and notation}\label{s2}


Fix $\Delta\notin\R$ and let $E=\R\cup\{\Delta\}$.
Define a metric $d$ on $E$ by
\[
d(x,y)=\cases{|x-y|\wedge1, &\quad $x,y\in\Bbb R$,\cr
1, &\quad $x\in\Bbb R, y=\Delta$,\cr
0, &\quad $x=y=\Delta$.}
\]
Thus $\Delta$ is an isolated point, which will act as a cemetery
state, and for $x, y\in\R$, $|x-y|\to0$ if and only if $d(x,y)\to0$.
Let $D$ be the Skorohod space of functions on $[0,\infty)$, taking
values in the metric space $E$, and which are
right continuous with left limits.
It is often convenient to assume that $X$ is given as the coordinate
process on $D$. We will interchangeably write
$X$ or $w$ depending on which seems clearer in the context. The usual
right continuous completion of the filtration generated by the
coordinate maps will be denoted by $\{\mathcal{F}_t\}_{t\ge0}$.
$P_z$ denotes the probability measure induced on
$\mathcal{F}=\bigvee_{t\ge0} \mathcal{F}_t$
by the L\'evy process\vspace*{1pt} starting at $z\in\R$.
We sometimes write just $P$ for $P_0$.
The shift operators $\theta_t\dvtx D\to D$, $t\ge0$,
are defined by $(\theta_t(w))_s=w(t+s)$.

For a given function
$w=(w_t)_{t\ge0}\in D$,
and
$r\ge0$, let
$w_{[0,r)}=(w_{[0,r)}(t))_{t\ge0}\in D$ denote the killed path
\[
w_{[0,r)} (t)=\cases{w_t, &\quad $0\le t<r$,\cr
\Delta, &\quad $t\ge r$.}
\]
For any $\rho\dvtx D\to[0,\infty]$ we then have the corresponding element
$\wta\in D$ defined by $\wta=w_{[0,\rho(w))}$.
For $x\in E$, let $c^x \in D$ be the constant path $c^x_t=x$ for all
${t\ge0}$.
If $w, w'\in D$, $w-w'$ denotes the path in $D$ given by
\[
(w-w')_t=
\cases{
w_t-w_t', &\quad if $t<\tau_\Delta(w)\wedge\tau_\Delta(w')$,\cr
\Delta, &\quad otherwise.}
\]
Let
\[
\tau_z=\tau_z(w)=\inf\{t>0\dvtx w_t>z\},\qquad \tau_\Delta= \tau
_\Delta(w) = \inf\{t>0\dvtx w_t=\Delta\}.
\]
For notational convenience, we will interchangeably write $w_t$ and
$w(t)$, $\tau_z$ and $\tau(z)$, etc.
Observe that for any $t\ge0$ and $w\in D$
%
%
\begin{equation}\label{kps}
\tau_\Delta\bigl(\wt\bigr)=t \qquad\mbox{if }\tau_\Delta(w)\ge t.
\end{equation}

We adopt the following
notation from \cite{GM1} which is very standard in the area;
cf. \cite{bert,doneystf} and \cite{kypbook}.
Let $(L_s)_{s\ge0}$ denote the local time of $X$ at its maximum, and
$(L^{-1}_s,H_s)_{s \geq0}$
the weakly
ascending bivariate ladder
process.
When $X_t\to-\infty$ a.s., $L_\infty$ has an exponential\vadjust{\goodbreak}
distribution with some parameter $q>0$, and the defective process
$(L^{-1},H)$ may
be obtained from a nondefective process
$(\mathcal{L}^{-1},\mathcal{H})$ by independent exponential killing at
rate $q > 0$.
Thus
%
%
\begin{equation}\label{H}
\bigl((L^{-1}_s,H_s)\dvtx s<L_\infty\bigr)
\eqdr\bigl((\mathcal{L}_s^{-1},\mathcal{H}_s)\dvtx s<e(q)\bigr),
\end{equation}
where $e(q)$ is independent of $(\mathcal{L}^{-1},\mathcal{H})$ and
has exponential
distribution with parameter $q$.

We denote the bivariate L\'{e}vy measure of
$(\mathcal{L}^{-1},\mathcal{H})$
by $\Pi_{{L}^{-1},{H}}(\cdot,\cdot)$.
The Laplace exponent $\kappa(a,b)$
of $(L^{-1},H)$,
defined by
%
%
\begin{equation} \label{kapdef}
e^{-\kappa(a,b)} = E (e^{-a{L}^{-1}_1 -b{H}_1}; L_\infty>1) = e^{-q}E
e^{-a\mathcal{L}^{-1}_1 -b\mathcal{H}_1}
\end{equation}
for values of $a,b\in\R$ for which the expectation is finite,
may be written
%
%
\begin{equation} \label{kapexp}\qquad
\kappa(a,b)
=
q+\rmd_{L^{-1}}a+\rmd_Hb+\int_{t\ge0}\int_{x\ge0}
(1-e^{-at-bx})
\Pi_{{L}^{-1}, {H}}(\rmd t, \rmd x),
\end{equation}
where $\rmd_{L^{-1}}\ge0$ and $\rmd_H\ge0$ are drift constants.
The bivariate renewal function of
$(L^{-1},H)$ is
%
%
\begin{equation}\label{Vkdef}
V(t,x)= \int_0^\infty e^{-qs}P(\mathcal{L}_s^{-1}\le t,\mathcal
{H}_s\le x)\,\rmd s.
\end{equation}
Its Laplace transform is given by
%
%
\begin{equation}\label{Vkap}\qquad
\int_{t\ge0}\int_{x\ge0}e^{-at-bx} V(\rmd t,\rmd x)=\int_{s\ge0}
e^{-qs}E(e^{-a \mathcal{L}_s^{-1}-b\mathcal{H}_s})\,\rmd s =\frac
{1}{\kappa(a,b)},
\end{equation}
provided $\gk(a,b)>0$. We will also frequently consider the renewal
function of $H$, defined on $\R$
by
%
%
\begin{equation}\label{VHdef}
V(x)= \int_0^\infty e^{-qs} P(\mathcal{H}_s\le x)\,\rmd s
=\lim_{t\to\infty}V(t,x).
\end{equation}
Observe that $V(x)=0$ for $x<0$, while $V(0)>0$ iff $\mathcal{H}$ is
compound Poisson. Also
%
%
\begin{equation}\label{Vinf}
V(\infty):=\lim_{x\to\infty} V(x)=q^{-1}.
\end{equation}

Let $\whX_t=-X_t$, $t\ge0$, denote the dual process and $(\whL^{-1},
\whH)$ the corresponding \textit{strictly} ascending bivariate
ladder processes of $\whX$. This is the same as the weakly ascending
process if $\whX$ is not compound Poisson.
All quantities relating to $\whX$ will be denoted in the obvious way,
for example,
$\Pi_{{\whL}^{-1}, {\whH}}$, $\widehat\kappa$ and $\whV$.
The reason for this choice of $(\whL^{-1}, \whH)$ is that we may
then, for any L\'evy process,
choose the normalization of the local times $L$ and $\whL$ so that
the Wiener--Hopf factorization takes the form
%
%
\begin{equation}\label{WH}
\kappa(a,-ib) \wh\kappa(a,ib)
= a-\Psi_X( b),\qquad a\ge0, b\in\R.
\end{equation}

Throughout the paper our principal assumption will be
(\ref{c1}). In that case, by Proposition 5.1 of \cite{kkm},
%
%
\begin{equation}\label{kpos}
\gk(a,-\ga)>0 \qquad\mbox{for $a\ge0$.}\vadjust{\goodbreak}
\end{equation}
Furthermore, by analytic extension, it follows from (\ref{WH}) that
%
%
\begin{equation}\label{WH1}
{\kappa(a,-z) \wh\kappa(a,z)
= a-\Psi_X(-\rmi z) \qquad\mbox{for } a\ge0, 0\le\Re z\le\ga.}
\end{equation}
Set
%
%
\begin{eqnarray}\label{const}
\ca&=& -\ln Ee^{\alpha X_1}=
-\Psi_X(-\rmi\alpha)=
\kappa(0,-\alpha)\whk(0,\alpha),\nonumber\\[-8pt]\\[-8pt]
\cb&=&\frac{\kappa(0,-\alpha)}{q},\qquad
\cab=\ca\cb.\nonumber
\end{eqnarray}
Note that $\ca, \cb\in(0,\infty)$ under (\ref{c1}).
These constants appear in several formulas throughout the paper.
For future reference we also note that
%
%
\begin{equation}\label{c1=}
\ca\int_0^\infty Ee^{\alpha X_t} \,\rmd t = \ca\int_0^\infty
(Ee^{\alpha X_1})^t \,\rmd t = 1
\end{equation}
and, letting $\Vbar(z)=V(\infty)-V(z)$, $z\in\R$, we have by (\ref
{Vkap}) and (\ref{Vinf})
%
%
\begin{eqnarray}\label{c2=}
\cb\int_z \alpha e^{-\alpha z} q\Vbar(-z) \,\rmd z
&=&\cb\biggl(1 +
\int_{z\ge0} \alpha e^{\alpha z}q\Vbar(z) \,\rmd z\biggr)\nonumber\\
&=& \cb q\int_{z\ge0} e^{\alpha z}V(\rmd z)\\
&=& \frac{\cb q}{\gk(0,-\ga)}=1.\nonumber
\end{eqnarray}

The following important asymptotic estimate can be found in \cite{kkm}.
Assuming~(\ref{c1}),
%
%
\begin{equation}\label{clim}
\lim_{u\to\infty}
\frac{\pibar_X^{+}(u)}{q\Vbar(u)} =\cab.
\end{equation}
This provides information about the probability of eventual ruin through
the \textit{Pollacek--Khintchine formula}
%
%
\begin{equation}\label{PK}
P\bigl(\tau(u)<\infty\bigr)=q\Vbar(u).
\end{equation}
A further useful estimate from \cite{kkm}, holding under (\ref{c1}),
is
%
%
\begin{equation}\label{XHlim}
\lim_{u\to\infty}
\frac{\pibar_X^{+}(u)}{\pibar_H(u)} =\whk(0,\alpha)\in(0,\infty),
\end{equation}
where $\Pi_H$ is the L\'evy measure of $H$, and $\pibar_H$ is its tail.
In particular, this implies
%
%
\begin{equation}\label{HSga}
\overline{\Pi}_{H} \in\mathcal{S}^{(\alpha)},
\end{equation}
since $\mathcal{S}^{(\alpha)}$ is closed under tail equivalence
(see \cite{EG}, Theorem 2.7).


\section{Main results}\label{s3}

We next introduce the basic components of the limiting process,
namely, processes $W$ and $Z$, and a random variable $\rho$.
These\vadjust{\goodbreak}
three random elements are independent.
The distribution of $W$ is given by
%
%
\begin{equation}\label{defW}
P(W\in\rmd w) = \cb\int_{z\in\R}\alpha e^{-\alpha z} q\Vbar
(-z)\,\rmd z
P_z\bigl(X\in\rmd w|\tau(0) <\infty\bigr)
\end{equation}
for $w\in D$.

\noindent[Recall that $\Vbar(y)=q^{-1}$ for $y<0$.]
Thus $W$ has the law of $X$ conditioned on $\tau(0)< \infty$
and started with initial distribution
%
%
\begin{equation}\label{W0}
P(W_0\in\rmd z)=\cb\alpha e^{-\alpha z} q\Vbar(-z) \,\rmd z,\qquad z\in\R.
\end{equation}
Observe that (\ref{W0}) is indeed
a probability distribution by (\ref{c2=}).
Let $Z$ be the Esscher transform of $X$, defined by
%
%
\begin{eqnarray}\label{Z}
&&
P(\{Z_t\dvtx 0\le t\le s\}\in B_s, Z_s\in\rmd x)\nonumber\\
&&\qquad=\frac{e^{\alpha x}P(\{
X_t\dvtx 0\le t\le s\}\in B_s, X_s\in\rmd x)}{Ee^{\alpha X_s}}\\
&&\qquad=e^{\ca s}{e^{\alpha x}P(\{X_t\dvtx 0\le t\le s\}\in B_s, X_s\in\rmd x)},
\nonumber
\end{eqnarray}
where $B_s$ is a Borel set in $\R^{[0,s]}$ and $x\in\R$.
Finally let
%
%
\begin{equation}\label{texp}
\rho\mbox{ be exponentially distributed with parameter } \ca.
\end{equation}
%

Let $H\dvtx D\otimes D\to\R$ be measurable with respect to the product
$\sigma$-algebra and set
%
%
\begin{equation}\label{Gdef}
G(w,z)= E_z[H(w,X);\tau(0) <\infty],\qquad w\in D, z\in\R.
\end{equation}
We denote by $\mathcal{H}$ the class of such functions $H$ which satisfy
%
%
\begin{eqnarray} \label{f1}
&H(w,w')e^{\gt w_{\tD-}I(w_{\tD-}\le0)} \mbox{ is bounded
for some
$\gt\in[0,\ga)$};&
\\
%
%
\label{f2}
&G(w,\cdot) \mbox{ is
continuous a.e. on $ (-\infty,\infty)$ for every $w\in D$}.&
\end{eqnarray}
%
For example, if $H$ is bounded and continuous in the product Skorohod
topology on
$D\otimes D$, these conditions hold with $\gt=0$. More general
conditions on $H$, which
ensure that (\ref{f2}) holds, will be discussed below.
Taking $\gt>0$ in (\ref{f1}) allows for certain unbounded functions
$H$, which will be used in Section \ref{s8}.

Here is our main theorem.
%
%
\begin{theorem}[(Path decomposition)]\label{THML1}
$\!\!\!$Assume (\ref{c1}). Then for any \mbox{$H\!\in\!\mathcal{H}$}
%
%
\begin{equation}\label{L1}\qquad
\lim_{x\to\infty}\lim_{u\to\infty}
E^{(u)}H\bigl(\X,\XD-c^u\bigr)=EH\bigl(Z_{[0, \rho)} ,W\bigr).
\end{equation}
\end{theorem}
%

The reason for introducing $x$ and taking the limit is to capture the
difference in behavior of the conditioned process before and after
entering a~neighborhood of $u$.
The heuristic meaning of the result is that
the conditioned process, for large $u$, can be approximated as follows:
\begin{itemize}
\item run the process $Z$ for times $0\le t<\rho$;
\item then, run the process $u+W$ from time $\rho$ on,
that is, at time $\rho+t$ the value of the process is $u+W_t$.\vadjust{\goodbreak}
\end{itemize}
Thus the process behaves like $Z$ up until an independent exponential
time~$\rho$, at which time it makes a large jump from a neighborhood
of $0$ to a~neighborhood of~$u$. Its position immediately prior to the
jump is $Z_{\rho-}$ and its position after the jump is $u+W_0$. If
$W_0>0$, the process $X-u$ behaves like~$X$ started at $W_0$. If
$W_0\le0$, the process $X-u$ behaves like $X$ started at $W_0$ and
conditioned on $\tau(0)<\infty$.
This behavior is significantly different from the Cram\'er and
subexponential cases discussed earlier.

It is apparent that many asymptotic results will
flow from Theorem \ref{THML1}.
We develop some of these in Sections \ref{s5}--\ref{s8}.
The literature to date has focused on deficit at ruin (overshoot)
and surplus prior to ruin (undershoot).
We use Theorem \ref{THML1} to derive these and related results in
Section \ref{s7}.
Of perhaps greater importance in insurance risk theory,
though, is the probability of ruin occurring
in finite time.
So far this has been neglected in studies of this type
(except, see the paper of Braverman \cite{br} discussed below).
%
The following result, derived from Theorem \ref{THML1}, gives a
completely explicit representation of the asymptotic distribution of
the ruin time.
%
%
\begin{theorem}[(Asymptotic distribution of ruin time)]\label{cor1}
Assume (\ref{c1}). Then for $t\ge0$
%
%
\begin{eqnarray}\label{RT}
\lim_{u\to\infty}
P^{(u)}\bigl(\tu\le t\bigr)&=&P\bigl(\rho+\tW(0)\le
t\bigr)\nonumber\\[-8pt]\\[-8pt]
&=&\cb E (e^{\alpha\Xbar_{t-\rho}}; \rho\le t),
\nonumber
\end{eqnarray}
where $\rho$ is independent of $X$ and $W$ and has exponential
distribution with parameter
$\ca$.
\end{theorem}

We can compare this result with those of Braverman \cite{br}.
He assumes, as we do, that
$\pibar_{X}^+(x)\in\mathcal{S}^{(\alpha)}$ for an $\alpha> 0$,
and his Theorem 2.1 can be used to deduce that
\[
\lim_{u\to\infty} \frac{P(\tu\le t)}{\pibar_{X}^+(u)}
\]
exists for each $t>0$, and hence, via (\ref{clim}), that
$\lim_{u\to\infty} P^{(u)}(\tu\le t)$
also exists. However, the expressions
thus obtained for these limits
are highly inexplicit, and it is not at all clear from them
whether or not the limiting distribution is proper (total mass 1).
Theorem \ref{cor1} gives
a much simpler expression for the limiting distribution
and establishes that it is indeed proper,
being the convolution of two proper probability distributions.


\section{Proof of path decomposition}\label{s4}

Let $\mathcal{B}$ denote the Borel sets on $\R$, $\mathcal
{B}([0,\infty))$ the
Borel sets on $[0,\infty)$ and set $\mathcal{D}= D\otimes[0,\infty
) \otimes(-\infty,\infty)$.
For~$K\in(-\infty,\infty]$ and $x\in[0,\infty]$, define
measures
$\mu_{K}$ and $\nu_x$ on $\mathcal{F}\otimes\break \mathcal{B}([0,\infty
))\otimes\mathcal{B}$ by
%
%
\begin{equation}\label{muK}\quad
\mu_K(\rmd w, \rmd t, \rmd\phi) = \ca I(\phi< K)e^{\alpha\phi
}P\bigl(\Xt\in\rmd w; X_{{t}-}\in\rmd\phi\bigr) \,\rmd t\vadjust{\goodbreak}
\end{equation}
and
%
%
\begin{equation}\label{nuK}\quad
\nu_x(\rmd w', \rmd r, \rmd z) = \cb I(z>-x)\alpha e^{-\alpha z}\,\rmd
z P_z\bigl(X\in\rmd w';\tau(0) \in\rmd r \bigr).
\end{equation}
We will write $\mu$ and $\nu$ for $\mu_\infty$ and $\nu_\infty$,
respectively.
These are finite measures and indeed $\mu$ and $\nu$ are probability
measures on $\mathcal{D}$ since by (\ref{c1=})
%
%
\begin{equation}\label{mkt}
\mu(\mathcal{D})=\ca\int_0^\infty Ee^{\alpha X_t} \,\rmd t=1,
\end{equation}
and by (\ref{c2=}) and (\ref{PK})
%
%
\begin{eqnarray}\label{nkt}
\nu(\mathcal{D})&=&\cb\biggl(1+\int_{z\le0} \alpha e^{-\alpha z}
P_z\bigl(\tau(0)< \infty\bigr) \,\rmd z\biggr)\nonumber\\
&=&\cb\biggl(1+\int_{z\ge0} \alpha e^{\alpha z} P\bigl(\tau(z)< \infty\bigr)
\,\rmd z\biggr)\\
&=&\cb\biggl(1 + \int_{z\ge0} \alpha e^{\alpha z}q\Vbar(z) \,\rmd
z\biggr)=1.\nonumber
\end{eqnarray}
In a slight abuse of notation we will denote the marginal measures in
the obvious way. Thus, for example,
%
%
\begin{eqnarray}\label{munu}\quad
\mu_K(\rmd w,\rmd\phi)&=&\ca\int_0^\infty I(\phi< K)e^{\alpha\phi
}P\bigl(\Xt\in\rmd w; X_{{t}-}\in\rmd\phi\bigr) \,\rmd
t,\nonumber\\[-8pt]\\[-8pt]\quad
\nu_x(\rmd w') &=& \cb\int_{z>-x}\alpha e^{-\alpha z}\,\rmd z P_z\bigl(X\in
\rmd w';\tau(0) <\infty
\bigr).\nonumber
\end{eqnarray}
%

From (\ref{mkt}) and (\ref{nkt}), $\mu(\rmd w)$ and $\nu(\rmd w')$
define probability measures on $D$.
From (\ref{PK}) and (\ref{defW}), it is clear that
$\nu(\rmd w')=P(W\in\rmd w')$.
The following result identifies $\mu$
as the distribution of
$Z_{[0, \rho)}$, where $Z$ and $\rho$ are given by~(\ref{Z}) and (\ref{texp}), respectively.
%
%
\begin{prop}\label{PYZ}
Let $\tZ$ have law given by $P(\tZ\in\rmd w)=\mu(\rmd w)$, and set
$\tau_{\tZ}=\tau_{\Delta}(\tZ)=\inf\{t>0\dvtx\tZ_t=\Delta\}$. Then
with $\rho$ and $Z$ as above,
%
%
\begin{equation}\label{YZ}
\{\tZ_t\dvtx t<\tau_{\tZ}\}\stackrel{d}{=}\{Z_t\dvtx t<\rho\}.
\end{equation}
\end{prop}
\begin{pf}
For any $B_s\in\mathcal{B}([0,s])$
\begin{eqnarray*}
&&P(\{\tZ_t\dvtx 0\le t\le s\}\in B_s, \tZ_s\in\rmd x, s<\tau_{\tZ})\\
&&\qquad=\ca\int_{r>s}\int_\phi e^{\alpha\phi}P(\{X_t\dvtx 0\le t\le s\}\in
B_s, X_s\in\rmd x, X_{{r}-}\in\rmd\phi) \,\rmd r\\
&&\qquad=\ca P(\{X_t\dvtx 0\le t\le s\}\in B_s, X_s\in\rmd x)\int_{r>0}\int_\phi
e^{\alpha\phi}P(X_{{r}-}\in\rmd\phi- x) \,\rmd r\\
&&\qquad=e^{\alpha x}P(\{X_t\dvtx 0\le t\le s\}\in B_s, X_s\in\rmd x)
\qquad\mbox{[by (\ref{c1=})]}\\
&&\qquad=P(\{Z_t\dvtx 0\le t\le s\}\in B_s, Z_s\in\rmd x){e^{-\ca s}}
\qquad\mbox{[by (\ref{Z})]}\\
&&\qquad=P(\{Z_t\dvtx 0\le t\le s\}\in B_s, Z_s\in\rmd x, s<\rho).
\end{eqnarray*}
Integrating out $x$ completes the proof.
\end{pf}
%
%
\begin{lemma}\label{L2}
Fix $x \in[0, \infty)$, $u>x$,
$A\subset(-\infty,u-x]$ and $B\subset(u-x,\infty)$.
Then for any $H\in\mathcal{H}$ which is nonnegative,
%
%
\begin{eqnarray}\label{L2a}\quad
&&E\bigl[H\bigl(\X,\XD-c^u\bigr)\dvtx X_{{{\tux}-}}\in A, X_{{{\tux}}}\in B, \tu<\infty\bigr]
\nonumber\\
\quad&&\qquad=\int_0^\infty\rmd t\int_{w\in{D}}\int_{\phi\in A}\int_{{z\in B-u}}
G(w, z)\Pi_X^+(u-\phi+\rmd z)\\
\quad&&\hspace*{145.7pt}{}\times P\bigl(\Xt\in\rmd w, X_{t-}\in\rmd\phi,
\tux\ge t\bigr),
\nonumber
\end{eqnarray}
where $G$ is defined by (\ref{Gdef}).
\end{lemma}
\begin{pf}
By the strong Markov property
\begin{eqnarray*}
&&E\bigl[H\bigl(\X,\XD-c^u\bigr)\dvtx X_{{{\tux}-}}\in A, X_{{{\tux}}}\in B, \tu<\infty\bigr]
\\
&&\qquad=
E\bigl[G\bigl(\X, X_{{\tux}}-u\bigr)\dvtx X_{{{\tux}-}}\in A,X_{{{\tux}}}\in
B,\\
&&\hspace*{238.6pt}\tux<\infty\bigr].
\end{eqnarray*}
Since $AB=\varnothing$, $\Delta X_{{{\tux}}}>0$ on
$\{X_{{{\tux}-}}\in A, X_{{{\tux}}}\in B\}$.
Thus by the compensation formula
(see \cite{bert}, page 7)
\begin{eqnarray*}
&&E\bigl[G\bigl(\X, X_{{\tux}}-u\bigr)\dvtx X_{{{\tux}-}}\in A, X_{{{\tux}}}\in B, \tux
<\infty\bigr]\\
&&\qquad=E\sum_t G\bigl(\Xt, X_{{t}-}+\Delta X_t-u\bigr)I\bigl(X_{t-}\in A, \tux
\ge t\bigr)\\
&&\hspace*{56pt}{}\times I(X_{t-}+\Delta X_t\in B)\\
&&\qquad=E\int_0^\infty\rmd t\int_{\xi} G\bigl(\Xt, X_{{t}-}+\xi-u\bigr)
I\bigl(X_{t-}\in
A,\tux
\ge t\bigr)\\
&&\hspace*{84pt}{}\times I(X_{{t}-}+\xi\in B)\Pi_X^+(\rmd\xi)\\
&&\qquad=\int_0^\infty\rmd t\int_{w\in{D}}\int_{\phi\in A}\int_{\xi
+\phi\in B}
G(w,\phi+\xi-u)\Pi_X^+(\rmd\xi)\\
&&\hspace*{147pt}{}\times P\bigl(\Xt\in\rmd w,X_{t-}\in\rmd
\phi,\tux\ge t\bigr),
\end{eqnarray*}
and this is (\ref{L2a}).
\end{pf}

In conjunction with Lemma \ref{L2}, it is useful to note that,
for $u>x\ge0$,
%
%
\begin{equation}\label{nc}
P\bigl(X_{{{\tux}-}}<u-x, X_{{{\tux}}}=u-x, \tux<\infty\bigr)=0
\end{equation}
(see, e.g.,
\cite{GM1}, Lemma 5.1).\vadjust{\goodbreak}

We need two further observations before we come to the proof of Theorem~\ref{THML1}.
From (\ref{f1}) and (\ref{f2}), it follows immediately that for any $K$
%
%
\begin{equation}\label{f2a}
G(w,z)e^{\gt w_{\tD-}I(w_{\tD-}<K)}
\end{equation}
is bounded as a function of $(w,z)$ and
continuous a.e. in $z
$ for every $w\in D$.

Referring to (\ref{Salph}),
an important global bound on convolution equivalent distributions is
obtained by applying Theorem 1.5.6(ii) of \cite{BGT} to the function
\[
l(r)=(r\vee e)^{\ga} \overline{F}\bigl(\ln(r\vee e)\bigr),
\]
which is slowly varying as $r\to\infty$.
This yields the following version of Potter's bounds for regularly
varying functions. Assume (\ref{Salph}); then
for every $\veps>0$ there exists an $A=A_{\veps}$ such that
%
%
\begin{equation}\label{Pot}\qquad
\frac{\overline{F}(u+x)}{\overline{F}(u)}\le
A\bigl[e^{-(\ga-\veps)x}\vee e^{-(\ga+\veps)x}\bigr]
\qquad\mbox{for all } u\ge1, x\ge1-u.
\end{equation}
Clearly this estimate also applies to $\pibar_X^+$ since we may take
$\overline{F}(x)=\pibar_X^+(x)/\allowbreak\pibar_X^+(1)$ for $x\ge1$.
Similarly for
$\pibar_H$, since recall $\overline{\Pi}_{H}
\in\mathcal{S}^{(\alpha)}$ from (\ref{HSga}).

The key step in the proof of Theorem \ref{THML1} is the following
result.
%
%
\begin{theorem}\label{THM1}
Assume (\ref{c1}), and
fix $x \in[0, \infty)$ and $K\in(-\infty,\infty)$.
Then for any $H\in\mathcal{H}$
%
%
\begin{eqnarray}\label{1}
&&\lim_{u\to\infty}E^{(u)}\bigl[H\bigl(\X,\XD-c^u\bigr);X_{{{\tux}-}}<
K\bigr]\nonumber\\[-8pt]\\[-8pt]
&&\qquad=\int_{w\in{D}}\int_{w'\in{D}} H(w,w')\mu_K(\rmd w)\nu_x(\rmd
w').\nonumber
\end{eqnarray}
\end{theorem}
\begin{pf}
We first show that the expression for the limit is finite. By
Proposition~\ref{PYZ}, and independence of $Z_{[0, \rho)}$ and $W$,
%
%
\begin{equation}\label{ZtW}
P\bigl(Z_{[0, \rho)}\in\rmd w,W\in\rmd w'\bigr)=\mu(\rmd
w)\otimes\nu(\rmd w').
\end{equation}
Hence, using (\ref{munu}),
%
%
\begin{eqnarray}\label{fi}
&&
\int_{w\in{D}}\int_{w'\in{D}}|H(w,w')|\mu_K(\rmd w)\nu_x(\rmd w')\nonumber\\
&&\qquad\le\int_{w\in{D}}\int_{w'\in{D}}|H(w,w')| \mu(\rmd w)\nu(\rmd
w')\\
&&\qquad=E\bigl|H\bigl(Z_{[0, \rho)} ,W\bigr)\bigr|.\nonumber
\end{eqnarray}
To verify that the final expectation is finite, it suffices by
(\ref{f1}) to show that $Ee^{-\gt Z_{\rho-}I(Z_{\rho-}\le0)}<\infty
$. But by
(\ref{Z})
\[
Ee^{-\gt Z_{\rho-}}=\int_0^\infty Ee^{-\gt Z_{s}}P(\rho\in\rmd
s)=\ca\int_0^\infty Ee^{(\ga-\gt) X_{s}}\,\rmd s <\infty,
\]
if $0\le\gt<\ga$.\vadjust{\goodbreak}

We now prove convergence. Take $u$ large enough that $K<u-x$, and set
$A=(-\infty,K)$ and $B=(u-x,\infty)$ in (\ref{L2a}).
Then, recalling (\ref{nc}), we have
\begin{eqnarray*}
&&E\bigl[H\bigl(\X,\XD-c^u\bigr)\dvtx X_{{{\tux}-}}<K, \tu<\infty\bigr]
\\
&&\qquad=\int_0^\infty\rmd t\int_{w\in{D}}
\int_{\phi<K}\int_{{z >-x}} G(w, z)\Pi_X^+(u-\phi+\rmd z)\\
&&\qquad\quad\hspace*{110pt}{}\times P\bigl(\Xt\in\rmd w, X_{t-}\in\rmd\phi, \tux\ge t\bigr).
\end{eqnarray*}
Using that $K$ and $ x$ are fixed, and that as $u\to\infty$, $I(\tux
\ge t)\to1$ and
%
%
\begin{equation}\label{wc}
\frac{\Pi_X^+(u-\phi+\rmd z)}{\pibar_X^+(u)} \to e^{\alpha\phi}
\alpha e^{-\alpha z} \,\rmd z
\end{equation}
in the sense of weak convergence on $(-x,\infty)$, we will show
%
%
\begin{eqnarray}\label{Glt}\qquad
&&\int_0^\infty\rmd t\int_{w\in{D}}\int_{\phi<K}\int_{{z >-x}}
G(w, z)\frac{\Pi_X^+(u-\phi+\rmd z)}{\pibar_X^+(u)}\nonumber\\
&&\qquad\quad\hspace*{77.5pt}{}\times P\bigl(\Xt\in\rmd w, X_{t-}\in\rmd\phi,\tux\ge t\bigr)\\
&&\qquad \to
\cab^{-1}\int_{w'\in{D}}\int_{w\in{D}}H(w,w') \mu_K(\rmd w)\nu
_x(\rmd
w').\nonumber
\end{eqnarray}
By (\ref{clim}) and (\ref{PK}), this
will complete the proof.

Let
\[
\Lambda_u(w, \phi) = \int_{z>-x}G(w,z)\frac{\Pi_X^+(u-\phi+\rmd
z)}{\pibar_X^+(u)}.
\]
For fixed $w$, $G(w,\cdot)$ is bounded and continuous a.e. by (\ref
{f2a}). Thus by (\ref{wc}), for fixed $(w,\phi)$
\[
\Lambda_u(w,\phi)\to\int_{z>-x}G(w,z)\alpha e^{\alpha(\phi-z)}
\,\rmd z =: \Lambda(w,\phi).
\]
Next write
\[
\frac{\Pi_X^+(u-\phi+\rmd z)}{\pibar_X^+(u)}
=\frac{\Pi_X^+(u-\phi+\rmd z)}{\pibar_X^+ (u-\phi-x)}\frac{\pibar
_X^+(u-\phi-x)}{\pibar_X^+ (u)}.
\]
The first term is a probability measure on $(-x,\infty)$.
For the second term,
choose $\theta$ to satisfy (\ref{f1}) and
fix $\veps>0$ so that $\gt+2\veps<\ga$. By
(\ref{Pot}), there exists an $A$ so that
%
%
\begin{equation}\label{Pibd}
\frac{\pibar_X^+(u-\phi-x)}{\pibar_X^+ (u)}\le A \bigl[e^{(\alpha-\veps
)(\phi+x)}\vee e^{(\alpha+\veps)(\phi+x)}\bigr],
\end{equation}
if $u\ge1$ and $\phi+x\le u-1$.
Now for $\phi< K$,
\begin{eqnarray*}
e^{(\alpha-\veps)\phi}\vee e^{(\alpha+\veps)\phi}
&\le& e^{(\alpha-\veps)\phi}I(\phi<0)+e^{(\alpha+\veps)K}e^{(\alpha
-\veps)\phi}I(0\le\phi<K)
\\
&\le& e^{(\alpha+\veps)|K|}e^{(\alpha-\veps)\phi}.
\end{eqnarray*}
Thus if $u_0=:(K+x+1)\vee1$, then by (\ref{f2a}),
for some constant $C$ depending on
$H, K$ and~$x$,
%
%
\begin{equation}\label{Lub}
{\sup_{u\ge u_0}}|\Lambda_u(w,\phi)|\le C e^{(\alpha-\veps)\phi}
e^{-\gt w_{\tD-}I(w_{\tD-}< K)}\quad
\mbox{all $w\in D$, $\phi< K$.}\hspace*{-45pt}
\end{equation}
In particular, since $\ga-\veps-\gt>\veps$, for every $t\ge0$
%
%
\begin{eqnarray}\label{Lub1}\quad
\sup_{u\ge u_0}\bigl|\Lambda_u\bigl(\Xt,X_{{t}-}\bigr)\bigr|I(X_{{t}-}<K)&\le& Ce^{(\ga
-\veps-\gt)X_{{t}-}}I(X_{{t}-}<K)\nonumber\\[-8pt]\\[-8pt]
\quad&\le& C_1 e^{\veps X_{{t}-}}I(X_{{t}-}<K),\nonumber
\end{eqnarray}
where $C_1=Ce^{(\ga-\veps-\gt)|K|}$.
Next observe that
\begin{eqnarray*}
\Phi_u(t)&=&\!\!:\int_{w\in{D}}\int_{\phi<K}\Lambda_u(w,\phi)P\bigl(\Xt
\in\rmd w, X_{{t}-}\in\rmd\phi;\tux
\ge t\bigr)\\
&=&E\bigl[\Lambda_u\bigl(\Xt,X_{{t}-}\bigr);X_{t-}<K, \tux
\ge t\bigr]\\
&\to& E\bigl[\Lambda\bigl(\Xt,X_{{t}-}\bigr);X_{t-}<K\bigr]=:\Phi(t)
\end{eqnarray*}
as $u\to\infty$, by bounded convergence using (\ref{Lub1}). Further,
again by (\ref{Lub1}), for any $t\ge0$
\[
{\sup_{u\ge u_0}}|\Phi_u(t)| \le C_1E[e^{\veps X_{{t}-}};X_{t-}<K]\le
C_1 (Ee^{\ga X_{{t}}})^{\veps/\ga}=C_1 (Ee^{\ga
X_{1}})^{\veps t/\ga},
\]
where recall $Ee^{\ga X_{1}}<1$ by (\ref{c1}). Thus dominated
convergence gives
%
%
\begin{equation}\label{Phl}
\int_0^\infty\Phi_u(t) \,\rmd t \to\int_0^\infty\Phi(t) \,\rmd t.
\end{equation}
This is equivalent to (\ref{Glt}) since the limit, which is expressed
in (\ref{Phl}) as an iterated integral, may be rewritten as
%
%
\begin{eqnarray}\label{rewr}\quad
\int_0^\infty\Phi(t) \,\rmd t&=& \int_t\rmd t\int_{w\in{D}}\int
_{\phi<K} \int_{z>-x}G(w,z)\alpha e^{\alpha(\phi-z)} \,\rmd z \nonumber\\
\quad&&\qquad\quad\hspace*{67.5pt}{}\times P\bigl(\Xt
\in\rmd w, X_{{t}-}\in\rmd\phi\bigr)\nonumber\\
\quad&=& \int_t\rmd t\int_{w\in{D}}\int_{\phi<K} \int
_{z>-x}E_{z}[H(w,X);\tau(0)<\infty]\alpha e^{\alpha(\phi-z)} \,\rmd
z\nonumber\\
\quad&&\qquad\quad\hspace*{67.5pt}{}\times  P\bigl(\Xt\in\rmd w, X_{{t}-}\in\rmd\phi\bigr)\\
\quad&=& \cb^{-1} \int_t\rmd t\int_{w\in{D}}\int_{\phi<K}e^{\alpha
\phi} P\bigl(\Xt\in\rmd w, X_{{t}-}\in\rmd\phi\bigr)\nonumber\\
\quad&&\qquad\quad\hspace*{32.6pt}{}\times \int_{w'\in{D}}
H(w,w')\nu_x(\rmd w')\nonumber\\
\quad&=& (\ca\cb)^{-1}\int_{w\in{D}}\int_{w'\in{D}} H(w,w') \mu
_K(\rmd w)\nu_x(\rmd
w').\nonumber
\end{eqnarray}
This calculation is justified by absolute convergence of the final
integral, proved earlier in (\ref{fi}).
\end{pf}
\begin{pf*}{Proof of Theorem \ref{THML1}}
Assume (\ref{c1}). Then using (\ref{munu}) and dominated convergence,
which is justified by (\ref{fi}), we have
\begin{eqnarray*}
&&\lim_{K,x\to\infty}\int_{w\in{D}}\int_{w'\in{D}}H(w,w')\mu
_K(\rmd w)\nu_x(\rmd w')\\
&&\qquad= \int_{w\in{D}}\int_{w'\in{D}}H(w,w') \mu(\rmd w)\nu(\rmd w').
\end{eqnarray*}
Thus by (\ref{1}) and (\ref{ZtW}),
%
%
\begin{eqnarray}\label{Hkx1}\qquad
&&\lim_{K,x\to\infty}\lim_{u\to\infty}E^{(u)}\bigl[H\bigl(\X,\XD
-c^u\bigr);X_{{{\tux}-}}<K\bigr]\nonumber\\[-8pt]\\[-8pt]
&&\qquad=EH\bigl(Z_{[0, \rho)} ,W\bigr).\nonumber
\end{eqnarray}
Taking $H\equiv1$ in (\ref{Hkx1}) gives
%
%
\begin{equation}
\lim_{K,x\to\infty}\lim_{u\to\infty}P^{(u)}\bigl(X_{{{\tux}-}}<K\bigr)=1.
\end{equation}
Since $H$ is bounded on $\{w_{\tD-}\ge K\}$ by (\ref{f1}), it follows that
%
%
\begin{eqnarray}\label{Hkx2}
&&\lim_{K,x\to\infty}\lim_{u\to\infty}E^{(u)}\bigl[H\bigl(\X,\XD
-c^u\bigr);\nonumber\\[-8pt]\\[-8pt]
&&\hspace*{165.5pt}X_{{{\tux}-}}\ge K\bigr]=0.\nonumber
\end{eqnarray}
Combining (\ref{Hkx1}) and (\ref{Hkx2}) then proves (\ref{L1}).
\end{pf*}
%
%
\begin{remark}
The limiting operations in this section are simpler than
those in
\cite{kkm}, not requiring the splitting of integrals over subdomains
and associated delicate estimations.
Further, many of the calculations do not require
the full force of the $\mathcal{S}^{(\alpha)}$ condition. In
particular, the proof of (\ref{Glt}) only uses
(\ref{Salph}) prior to equation (\ref{Phl}). At this point the
additional condition $E e^{\ga X_1}<1$ is needed to ensure
that dominated convergence applies. Thus the proof actually shows that
under (\ref{Salph}), if $H\in\mathcal{H}$ is such that (\ref{Phl})
holds, then for any $x\ge0$ and $K\in(-\infty, \infty)$
\begin{eqnarray*}
&&\lim_{u\to\infty}\frac{E[H(\X,\XD-c^u);X_{{{\tux}-}}< K, \tu
<\infty]}{\pibar_X^+(u)}
\\
&&
\qquad=\int_0^\infty\Phi(t)\,\rmd t.
\end{eqnarray*}
This is the case if, for example, $\Phi_u$ are dominated by an
integrable function on $[0,\infty)$.
If, in addition,
\[
\int_{w\in{D}}\int_{w'\in{D}}| H(w,w')|\tmu_K(\rmd w)\tnu_x(\rmd
w')<\infty,
\]
where
%
%
\begin{equation}\label{tmuK}\quad
\tmu_K(\rmd w, \rmd t, \rmd\phi) = I(\phi< K)e^{\alpha\phi}P\bigl(\Xt
\in\rmd w; X_{{t}-}\in\rmd\phi\bigr) \,\rmd t
\end{equation}
and
%
%
\begin{equation}\label{tnux}\quad
\tnu_x(\rmd w', \rmd r, \rmd z) = I(z>-x)\alpha e^{-\alpha z}\,\rmd z
P_z\bigl(X\in\rmd w';\tau(0) \in\rmd r \bigr),
\end{equation}
then the limit may be rewritten as
\[
\int_0^\infty\Phi(t)\,\rmd t=\int_{w\in{D}}\int_{w'\in{D}}
H(w,w')\tmu_K(\rmd w)\tnu_x(\rmd w')
\]
as demonstrated in (\ref{rewr}).
Comparing (\ref{tmuK}) and (\ref{tnux}) with (\ref{muK}) and (\ref
{nuK}), note that the constants $\ca$ and $\cb$ must be excluded since
they no longer need to be finite and nonzero.
\end{remark}

We briefly address conditions on $H$, beyond measurability,
which ensure that (\ref{f2}) holds.
It is natural that such
conditions should relate to some type of continuity of $H$.
We will assume that for each $w\in D$, $H(w,\cdot)$ is continuous from
below on $\{\tau_0(w')<\infty\}$ a.s. $P_z$ for every $z$, that is,
for all $w\in D, z\in\R$
%
%
\begin{equation}\label{ctsa}\quad
\lim_{\veps\downarrow0}H(w,w'-c^\veps)
= H(w,w') \qquad\mbox{a.s.
$P_z$ on $\{\tau_0(w')<\infty\}$}.
\end{equation}
%
This condition clearly holds if, for every $\omega\in D$, $H(w,\cdot
)$ is continuous in the uniform topology on $D$, and so, in particular,
if $H(w,\cdot)$ is continuous in any of the usual Skorohod topologies.
Several examples of functionals of interest that satisfy (\ref{ctsa})
are given in Lemma \ref{fcty}.
For a detailed discussion of the various topologies on Skorohod space,
see \cite{W}.
%
%
%
\begin{prop}\label{Hcond}
If $H$ is measurable and satisfies (\ref{f1}) and (\ref{ctsa}), then
(\ref{f2}) holds.
\end{prop}
\begin{pf}
For $y<z$, we have
\[
G(w,y)= E_y[H(w,X);\tau_0(X) <\infty]= E_z[H(w,X-c^{z-y});\tau
_{z-y}(X) <\infty].
\]
%
Next, by right continuity, $\tau_\veps(w')\downarrow\tau_0(w')
\mbox{ as } \veps\downarrow0$ for any $w'\in D$ with $\tau
_0(w')<\infty$.
Thus by (\ref{ctsa}), as $y\uparrow z$
%
%
\begin{eqnarray}\label{ctsaw}
&&H(w,X-c^{z-y})I\bigl(\tau_{z-y}(X)<\infty\bigr)\nonumber\\[-8pt]\\[-8pt]
&&\qquad\to H(w,X)
I\bigl(\tau_{0}(X)<\infty
\bigr)\qquad \mbox{a.s. } P_z.\nonumber
\end{eqnarray}
%
Hence, by bounded convergence, for each $w\in D$
\[
G(w,y)\to G(w,z)\vadjust{\goodbreak}
\]
%
as $y\uparrow z$. Thus $G(w,\cdot)$ is left continuous and
consequently continuous except at countably many points.
\end{pf}
%
%
\begin{remark}
Condition (\ref{ctsa}) can be weakened by requiring it to
hold except for a discrete set of $z$. This would result in $G(w,\cdot
)$ being left continuous except on a discrete set which again implies
continuity except at countably many points.
\end{remark}

One technical point should be mentioned.
In order that the expression in~(\ref{f2}) make sense,
$H(w,\cdot)$ must be measurable. If $H(w,\cdot)$ is continuous in the
uniform topology, this need not be the case since there are open sets
in the uniform topology which are not in the $\sigma$-algebra
generated by the coordinate maps
$\{w'_t\dvtx t\ge0\}$. This is why we impose the blanket condition that $H$
be measurable with respect to the product $\sigma$-algebra on
$D\otimes D$.

For later reference we note that $H(w,w')=e^{-\gt w_{\tD-}}$ trivially
satisfies~(\ref{ctsa}), and if
$\gt\in[0,\ga)$, then $H$ also satisfies~(\ref{f1}). Thus by
Proposition \ref{Hcond}, $H\in\mathcal H$. Hence, by taking
$x=K=0$ in Theorem \ref{THM1}, it follows that
%
%
\begin{equation}\label{undbd}
\limsup_{u\to\infty}\Eu e^{-\gt X_{\tu-}}<\infty
\qquad\mbox{for
every $\gt\in[0,\ga)$}.
\end{equation}
We will later show that the limit exists and evaluate it (see
Proposition \ref{GSeg0}).


\section{General marginal convergence results}\label{s5}

In this section we provide a~rec\-ipe for constructing conditional
limit theorems for the fluctuation variables, by
specializing Theorem \ref{THML1}.
This gives, in Theorem \ref{DK10}, joint convergence of the main
variables of interest in insurance risk.
Again we need some preliminary results.

By convention we set
$w'_{0-}=w'_0$ and $G_{0-}(w')=0$.
Also we define $\overline{w}'_t= \sup_{0\le s\le t}w'_s$.
%
%
\begin{lemma}\label{fcty}
Each of the following functions is continuous from below on $\{\tau
_0(w')<\infty\}$ a.s. $P_z$, for all $z$:
\[
w'_0, \tau_0(w'), G_{\tau(0)-}(w'), \overline{w}'_{\tau
(0)-},w'_{\tau(0)-}, w'_{\tau(0)}.
\]
\end{lemma}
\begin{pf}
Clearly $w'_0$ is continuous from below without any extra conditions.
Now assume that $\tau_0(w')<\infty$.
Let $\veps>0$ be sufficiently small that $\tau_\veps(w')<\infty$. Then
%
%
\begin{equation}\label{tcont2}
\tau_0(w'-c^\veps)=\tau_\veps(w').
\end{equation}
Thus by right continuity
%
%
\begin{equation}\label{tcont1}
\tau_0(w'-c^\veps)\downarrow\tau_0(w') \qquad\mbox{as } \veps
\downarrow0,
\end{equation}
which proves $\tau_0(w')$ is continuous from below on $\{\tau
_0(w')<\infty\}$.
Next, from~(\ref{tcont2}) we have
%
%
\begin{equation}\label{wv}
(w'-c^\veps)_{\tau_0(w'-c^\veps)\wedge\cdot}= w'_{\tau_{\veps
}(w')\wedge\cdot}-\veps.\vadjust{\goodbreak}
\end{equation}
If $w'_{\tau(0)}>0$, then $\tau_\veps(w')=\tau_0(w')$ for all
$\veps\in(0,w'_{\tau(0)})$, so the result for the remaining
functionals follows immediately from (\ref{wv}).
Thus we assume $w'_{\tau(0)}=0$, in which case $\tau_\veps(w')>\tau
_0(w')$ for all $\veps>0$. Now
$P_z(w'_{\tau(0)}=0)=0$ if $z>0$ so we only need to consider $z\le0$.
If $z<0$, then by Lemma 5.1 of~\cite{GM1}, $w'_{\tau_0(w')-}=w'_{\tau
_0(w')}$, and consequently also
$G_{\tau_0(w')-}(w')=\tau_0(w')$, a.s.~$P_z$.
This continues to hold for $z= 0$, since applying the strong Markov
property at time $\tau_0(w')$, shows $\tau_0(w')=0$ a.s. $P_0$ when
$w'_{\tau(0)}=0$.
Thus by right continuity, we have $P_z$ a.s.
\[
(w'-c^\veps)_{\tau_0(w'-c^\veps)-}=w'_{\tau_\veps(w')-}-\veps\to
w'_{\tau_0(w')}=w'_{\tau_0(w')-}
\]
and
\[
G_{\tau_0(w'-c^\veps)-}(w'-c^\veps)=G_{\tau_\veps
(w')-}(w')\downarrow G_{\tau_0(w')}(w')=\tau_0(w')=G_{\tau_0(w')-}(w').
\]
The proofs for the remaining functionals are similar.\vspace*{-2pt}
\end{pf}
%
%
\begin{remark}
The above result is false if we replace continuous from below with
continuous from above. For example,
if $X$ is a Poisson process, then for any $\veps>0$
\[
P_{0}\bigl(\tau_0(w'+c^\veps)=0\bigr)=1,\qquad P_{0}\bigl(\tau_0(w')=0\bigr)=0.\vspace*{-2pt}
\]
\end{remark}

It will be convenient to write
%
%
\begin{equation}\label{Ydef}
Y=\XD-c^u \qquad\mbox{if $\tux<\infty$.}
\end{equation}
Thus $Y_t= X_{t+\tau(u-x)}-u$, $t\ge0$, and, in particular, $Y_0=
X_{\tau(u-x)}-u$, when $\tux<\infty$.
Of course $Y=Y(u,x)$, but to simplify the notation, we suppress the
dependence on $u$ and $x$.
{}From Theorem \ref{THML1} we have that $Y$ converges to $W$ under
$\Pu$,
as $x,u\to\infty$, in the sense specified there.
Likewise, $X_{[0,\tau(u-x))}$ converges to $Z_{[0,\rho)}$
in the sense of Theorem \ref{THML1}, and in fact we have joint
convergence.
This provides us with a means for constructing limit theorems
for the fluctuation variables.
The first step is in the next proposition.
Recall the definition of $\Xbar$ in (\ref{GXbar}), and define $\Wbar
$ and $\Zbar$ analogously.
Note that in (\ref{Hfg1})
we replace the variables on the left-hand side with those
on the right-hand side in the limit, as just described. Since $Z$ is
a.s. continuous at $\rho$, one may further replace the subscripts
$\rho-$ by $\rho$ in (\ref{Hfg1}), but we leave them in their
present form to help emphasize the remark in the previous sentence.\vspace*{-2pt}
%
%
\begin{prop}\label{Pfg1}
Assume (\ref{c1}) and suppose
$f\dvtx\R^{10}\to\R$ is bounded, measurable and
jointly continuous in the last six arguments. Let $0\le\gt<\ga$ and set
\begin{eqnarray*}
H(w,w')&=&f\bigl(G_{\tD-}(w), \tau_\Delta(w), \overline{w}_{\tD-},
w_{\tD-}, w'_0,\\
&&\hspace*{12.4pt} G_{\tau(0)-}(w'), \tau_0(w'),
\overline{w}'_{\tau(0)-},w'_{\tau(0)-}, w'_{\tau(0)}\bigr)\\
&&{}
\times e^{-\gt w_{\tD-}I(w_{\tD-}\le0)}I\bigl(\tau
_\Delta(w)<\infty, \tau_0(w')<\infty\bigr).
\end{eqnarray*}
Then $H$ satisfies (\ref{ctsa}) and hence,
%
%
\begin{eqnarray}\label{Hfg1}
&&\lim_{x\to\infty}\lim_{u\to\infty}\Eu f\bigl(G_{\tux-}, \tux, \Xbar
_{\tux-},\nonumber\\[-2pt]
&&\hspace*{81pt} X_{\tux-},
Y_0, G^Y_{\tau(0)-}, \tau_0^Y, \Ybar_{\tau(0)-},Y_{\tau(0)-},
Y_{\tau(0)} \bigr)\nonumber\\[-2pt]
&&\quad\hspace*{35pt}{}
\times e^{-\gt X_{\tux-}I(X_{\tux-}\le0)}\\[-2pt]
&&\qquad =E f\bigl(G_{\rho-}^Z, \rho, \Zbar_{\rho-}, Z_{\rho-}, W_0,
G_{\tau(0)-}^W, \tau_0^W, \Wbar_{\tau(0)-},W_{\tau(0)-}, W_{\tau
(0)} \bigr)\nonumber\\[-2pt]
&&\qquad\quad\hspace*{0pt}{}\times e^{-\gt Z_{\rho-}I(Z_{\rho-}\le0)}.
\nonumber\vspace*{-3pt}
\end{eqnarray}
\end{prop}
\begin{pf}
$H$ satisfies (\ref{ctsa}) by Lemma \ref{fcty}. Thus by Proposition
\ref{Hcond} we may apply Theorem \ref{THML1}. Upon noting that $\tau
_\Delta(\X)=\tux$, the result then follows immediately.\vspace*{-3pt}
\end{pf}

In what is essentially a special case of the
description of the limiting process given after Theorem \ref{THML1},
we can immediately deduce from Proposition~\ref{Pfg1} that
the joint limiting distribution of the time of, the position prior to
and the position relative to $u$ after,
the large jump is that of $(\rho, Z_{\rho-}, W_0)$. To be precise,
under~(\ref{c1}), we have
%
%
\begin{eqnarray}\label{ex1}
&&\lim_{x\to\infty}\lim_{u\to\infty}P^{(u)}\bigl(\tux\in\rmd t,
X_{{{\tux}-}}\in\rmd\phi, X_{{\tux}}-u\in\rmd z\bigr)\nonumber\\[-2pt]
&&\qquad= P(\rho\in\rmd t, Z_{\rho-}\in\rmd\phi, W_0\in\rmd z)\\[-2pt]
&&\qquad=\ca e^{\alpha\phi}P(X_{{t}-}\in\rmd\phi) \,\rmd t
\cb\alpha e^{-\alpha z}q\Vbar(-z)\,\rmd z,
\nonumber
\end{eqnarray}
where the last equality follows from
(\ref{W0}), (\ref{Z}) and (\ref{texp}).
The exact meaning of this convergence is given by (\ref{L1}), which by
(\ref{Hfg1}), is stronger than the usual weak convergence.

Observe that on $(0,\infty)$,
$P(W_0\in\rmd z)=\cb\alpha e^{-\alpha z}\,\rmd z $ is the limiting distribution
of the overshoot $X_{\tu} -u$ when the overshoot is due to the
large jump. The limiting probability that the large jump results in an
overshoot of $u$ is $P(W_0>0)=\cb$.
A further discussion of the overshoot is given in Section \ref{s7}.
Note also that
(\ref{ex1}) describes the joint limiting distribution of the ruin
time, the claim surplus immediately prior to ruin and the shortfall at
ruin, when ruin is due to a large claim. This makes precise the
``intuitively obvious'' asymptotic independence observed after Theorem
11 in \cite{DK}, and extends it to also include the ruin time.

The next step in our recipe is to transfer from
the variables on the left-hand side of (\ref{Hfg1})
to the fluctuation variables. The key point is to
observe that if $\tu<\infty$ and $x<u$, then
%
%
\begin{eqnarray}\label{10c}
G_{\tu-}&=&G_{\tux-}I(Y_0>0)+\bigl({\tux} +G^Y_{\too-}\bigr)I(Y_0\le
0),\nonumber\hspace*{-15pt}\\[-2pt]
{\tu} -G_{\tu-}&=&\bigl({\tux} -G_{\tux-}\bigr)I(Y_0>0)\nonumber\hspace*{-15pt}\\[-2pt]
&&{}+\bigl({\tau_0^Y}
-G^Y_{\too-}\bigr)I(Y_0\le0),\nonumber\hspace*{-15pt}\\[-2pt]
X_{\tu}-u&=&Y_0I(Y_0>0)+Y_{\too}I(Y_0\le0)=Y_{\too},\nonumber\hspace*{-15pt}\\[-2pt]
\hspace*{34pt}\Xbar_{\tu-}-X_{\tu-}&=&\bigl(\Xbar_{\tux-}-X_{\tux-}\bigr)I(Y_0>0)\hspace*{-15pt}\\[-2pt]
&&{}+\bigl(\Ybar_{\too-}-Y_{\too-}\bigr)I(Y_0\le0),\nonumber\hspace*{-15pt}\\[-2pt]
\Xbar_{\tu-}&=&\Xbar_{\tux-}I(Y_0>0)+\bigl(u+\Ybar_{\too-}\bigr)I(Y_0\le 0),\nonumber\hspace*{-15pt}\\[-2pt]
u-\Xbar_{\tu-}&=&\bigl(u-\Xbar_{\tux-}\bigr)I(Y_0>0)-\Ybar_{\too-}I(Y_0\le0)\quad\mbox{and} \nonumber\hspace*{-15pt}\\[-2pt]
X_{\tu-}&=&X_{\tux-}I(Y_0>0)+\bigl(u+Y_{\too-}\bigr)I(Y_0\le0).\nonumber\hspace*{-15pt}
\end{eqnarray}
%

Since some limiting variables have mass at infinity, we will consider
weak convergence on $\R\cup\{\infty\}$. To be precise, we will
consider functions $f\dvtx\R^4\otimes(\R\cup\{\infty\})\to\R$ which
are jointly continuous in the sense that
$f(\mathbf{x}_n,{y}_n) \to f(\mathbf{x},{y})$ as $(\mathbf
{x}_n,{y}_n) \to(\mathbf{x},{y})$ for $\mathbf{x}_n,\mathbf{x}\in
\R^4$ and ${y}_n,{y}\in(\R\cup\{\infty\})$.
%
%
\begin{theorem}\label{DK10}
Assume (\ref{c1}).
Let $f\dvtx\R^4\otimes(\R\cup\{\infty\})\to\R$ be bounded and
jointly continuous. Then for $0\le\gt<\ga$
%
%
\begin{eqnarray}\label{101a}
&&\lim_{u\to\infty} \Eu f\bigl(G_{\tu-}, {\tu} -G_{\tu-},X_{\tu}-u,
\Xbar_{\tu-}-X_{\tu-}, \Xbar_{\tu-}\bigr)\nonumber\\
&&\qquad\hspace*{0pt}{}\times e^{-\gt X_{\tu-}I(X_{\tu
-}\le0)}\nonumber\\
&&\quad\qquad = E \bigl[f(G^Z_{\rho-},{\rho} -G^Z_{\rho-}, W_0, \Zbar_{\rho
-}-Z_{\rho-}, \Zbar_{\rho-})e^{-\gt
Z_{\rho-}I(Z_{\rho-}\le0)};
W_0> 0\bigr]\hspace*{-3pt}\\
&&
\qquad\qquad{} + E \bigl[f\bigl(\rho+G^W_{\too-},{\tau_0^W}
-G^W_{\too-},\nonumber\\
&&\quad\hspace*{70.5pt} W_{\tau
(0)}, \Wbar_{\tau(0)-}-W_{\tau(0)-}, \infty\bigr); W_0\le0\bigr]\nonumber
\end{eqnarray}
and
%
%
\begin{eqnarray}\label{101b}
&&\lim_{u\to\infty} \Eu f\bigl(G_{\tu-}, {\tu} -G_{\tu-},X_{\tu}-u,
\Xbar_{\tu-}-X_{\tu-}, u-\Xbar_{\tu-}\bigr)\nonumber\\
&&\qquad\hspace*{0pt}{}\times e^{-\gt X_{\tu-}I(X_{\tu
-}\le0)}\nonumber\\
&&\quad\qquad = E \bigl[f(G^Z_{\rho-},{\rho} -G^Z_{\rho-}, W_0, \Zbar_{\rho
-}-Z_{\rho-}, \infty)e^{-\gt Z_{\rho-}I(Z_{\rho-}\le0)}; W_0> 0\bigr]\\
&&\qquad\qquad{} + E \bigl[f\bigl(\rho+G^W_{\too-},{\tau_0^W} -G^W_{\too-},\nonumber\\
&&\hspace*{81pt} W_{\tau
(0)}, \Wbar_{\tau(0)-}-W_{\tau(0)-},- \Wbar_{\tau(0)-}\bigr);W_0\le
0\bigr].\nonumber
\end{eqnarray}
Thus, under $\Pu$,
%
%
\begin{eqnarray}\label{101aa}
&&\bigl(G_{\tu-}, {\tu} -G_{\tu-},X_{\tu}-u, \Xbar_{\tu-}-X_{\tu-},
\Xbar_{\tu-}\bigr)\nonumber\\
&&\qquad \to(G^Z_{\rho-},{\rho} -G^Z_{\rho-}, W_0, \Zbar_{\rho
-}-Z_{\rho-}, \Zbar_{\rho-})I(W_0> 0)\nonumber\\[-8pt]\\[-8pt]
&&\qquad\quad{} + \bigl(\rho+G^W_{\too-},{\tau_0^W} -G^W_{\too-},W_{\tau
(0)}, \Wbar_{\tau(0)-}-W_{\tau(0)-}, \delta_\infty\bigr)\nonumber\\
&&\hspace*{44pt}{}\times I(W_0\le
0)\nonumber
\end{eqnarray}
in the sense of weak convergence on $\R^4\otimes(\R\cup\{\infty\}
)$, and similarly for (\ref{101b}).
\end{theorem}
\begin{pf}
We only prove (\ref{101a}), as the proof of (\ref{101b}) is similar.
Write the expectation on the left-hand side of (\ref{101a}) as the sum
of two expectations, one over $Y_0>0 $ and the other over $Y_0\le0$.
Convergence of the expectation over $Y_0>0 $ to the first term on the
right-hand side of (\ref{101a}), as $u\to\infty$ then $x\to\infty
$, follows easily from (\ref{10c}) and Proposition \ref{Pfg1}, since
$P(W_0=0)=0$ and $f$ is bounded and jointly continuous. For the
expectation over $Y_0\le0$, first observe that if $Y_0\le0$, then
$X_{\tu-}=u+Y_{\too-}$, and so
on $\{Y_0\le0\}$
\[
e^{-\gt X_{\tu-}I(X_{\tu-}\le0)}=I\bigl(u+Y_{\too-}> 0\bigr)+ e^{-\gt X_{\tu
-}}I\bigl(u+Y_{\too-}\le0\bigr).
\]
Convergence of the expectation over $\{u+Y_{\too-}> 0, Y_0\le0\}$ to
the second term on the right-hand side of (\ref{101a}), as $u\to
\infty$ then $x\to\infty$, again follows from (\ref{10c}) and
Proposition \ref{Pfg1} since
%
%
\begin{equation}\label{uY0}
\lim_{x\to\infty}\lim_{u\to\infty}\Pu\bigl(u+Y_{\too-}\le0, Y_0\le
0 \bigr)=0.
\end{equation}
Since $f$ is bounded, it thus remains to show
\[
\lim_{x\to\infty}\lim_{u\to\infty}\Eu\bigl(e^{-\gt X_{\tu
-}};u+Y_{\too-}\le0, Y_0\le0 \bigr)=0.
\]
For this it suffices by (\ref{uY0}) and H\"older's inequality, to show
that for some $\gt'> \gt$,
\[
\limsup_{u\to\infty}\Eu e^{-\gt' X_{\tu-}}<\infty,
\]
which in turn holds for any $\gt'\in(\gt, \ga)$ by
(\ref{undbd}).\vspace*{-3pt}
\end{pf}

Theorem \ref{Pfg1} provides a general convergence
result for the variables of primary interest in insurance risk,
in the convolution equivalent case.
It contains and extends many previous results in the literature
as will be explained in Sections \ref{s6}--\ref{s8}.
The two components that make up the limiting distributions in Theorem
\ref{DK10} arise as a consequence of the process either overshooting
or undershooting the boundary
at the time of the large jump. We now give alternate characterizations
of these distributions in terms of quantities arising in fluctuation
theory.

Recall the definitions of $\gk$ and $V$ in (\ref{kapexp}) and (\ref
{Vkdef}), and of $\whk$ and $\whV$ in the paragraph following (\ref{Vinf}).
To avoid introducing further notation, there is clearly no harm in
assuming that the random elements $(W,Z,\rho)$ are independent of $X$.
In particular, $\rho$ has an exponential distribution with parameter
$\ca$ and is independent of $X$.
Then by the Wiener--Hopf factorization theorem, $(G_{\rho},\Xbar
_{\rho})$ and $(\rho-G_{\rho},\Xbar_{\rho}-X_{\rho})$ are
independent with
Laplace transforms given by
%
%
\begin{eqnarray}\label{WHfac}
Ee^{-aG_{\rho}-b\Xbar_{\rho}}&=&\frac{\gk(\ca,0)}{\gk(\ca
+a,b)},\nonumber\\[-9pt]\\[-9pt]
Ee^{-a(\rho-G_{\rho})-b(\Xbar_{\rho}-X_{\rho})}&=&\frac{\whk(\ca
,0)}{\whk(\ca+a,b)}\nonumber
\end{eqnarray}
for $a,b>0$ (see \cite{kypbook}, Section 6.4).\vadjust{\goodbreak}

Before stating the next result, we wish to make clear the meaning of
the notation $|V(\rmd t-r,z-\rmd y)|$ below. It is the measure defined
on Borel sets in~$\R^2$ by
\begin{eqnarray*}
&&\iint_{(t,y)} 1_A(t,y)|V(\rmd t-r,z-\rmd y)|\\[-2pt]
&&\qquad= \iint_{(t,y)}
1_A(t+r,z-y)V(\rmd t,\rmd y).
\end{eqnarray*}
Some authors omit the absolute values signs. We include them to
emphasize that the function $V(t-r,z-y)$ is increasing in $t$ and
decreasing in $y$, hence, the Stieltjes measure associated with it,
which assigns mass
\[
V(t_1-r,z-y_1)-V(t_1-r, z-y_0)-V(t_0-r, z-y_1)+V(t_0-r, z-y_0)
\]
to rectangles $(t_0,t_1]\times[y_0,y_1)$, is negative.
%
%
\begin{theorem}\label{DK11}
For $\gamma> 0, t\ge0,s\ge0, \theta\ge0, \phi\ge0$
%
%
\begin{eqnarray}\label{112}
&& P(G^Z_{\rho-}\in\rmd t, {\rho} -G^Z_{\rho-}\in\rmd s, W_0\in\rmd
\gamma, \Zbar_{\rho-}-Z_{\rho-}\in\rmd\phi,\nonumber\\
&&\hspace*{170pt} \Zbar_{\rho-}\in
\rmd\theta; W_0> 0)\\
&&\qquad = \cab\ga e^{-\ga(\gamma+\phi-\gt)} V(\rmd t,\rmd\theta
)\whV(\rmd s,
\rmd\phi)\,\rmd\gamma,\nonumber
\end{eqnarray}
where $\beta$ is given by (\ref{const}).

For $\gamma\ge0, t\ge0,s\ge0, v\ge0, y\ge0$
%
%
\begin{eqnarray}\label{102}
&& P\bigl(\rho+G^W_{\too-}\in\rmd t, {\tau_0^W} -G^W_{\too-}\in\rmd
s,\nonumber\\
&&\quad\hspace*{2pt}
W_{\tau(0)}\in\rmd\gamma, \Wbar_{\tau(0)-}-W_{\tau(0)-}\in\rmd v,
-\Wbar_{\tau(0)-}\in\rmd y; W_0\le0\bigr)\nonumber\\
&&\hspace*{2pt}\quad\qquad = \cab I(\gamma>0) \int_{r\le t} e^{-\ca r}\,\rmd r \int
_{z\ge y}\ga e^{\ga z}\,\rmd z |V(\rmd t-r,z-\rmd y)|\whV(\rmd s, \rmd
v)\\
&&\qquad\quad\hspace*{141pt}{}\times\Pi_X(\rmd\gamma+v+y)\nonumber\\
&&\hspace*{2pt}\qquad\qquad{} + \cab\,\rmd_H \int_{r\le t} e^{-\ca r}\,\rmd
r\int_{z\ge0}\ga e^{\ga z} V(\rmd t-r, \rmd z)\delta_0(\rmd s, \rmd
\gamma, \rmd v,\rmd
y),\nonumber
\end{eqnarray}
where $\delta_0$ denotes a point mass at the origin.
\end{theorem}
\begin{pf}
The form of the limit in (\ref{102}) follows from an extension of
Doney and Kyprianou's \cite{DK} quintuple law to include creeping, as
given in Griffin and Maller \cite{GM1}.
For $\gamma\ge0, t\ge0, r\ge0, v\ge0, y\ge0$, we have by (\ref
{defW}), and Theorems 3.1(ii) and 3.2 of \cite{GM1},
%
%
\begin{eqnarray}\label{10a}
&&P\bigl(G^W_{\too-}\in\rmd r, {\tau_0^W} -G^W_{\too-}\in\rmd s, W_{\tau
(0)}\in\rmd\gamma,\nonumber\hspace*{-15pt}\\
&&\hspace*{13.3pt}{} \Wbar_{\tau(0)-}-W_{\tau(0)-}\in\rmd v,
-\Wbar_{\tau(0)-}\in\rmd y; W_0\le0\bigr)\nonumber\hspace*{-15pt}\\
&&\qquad =\cb\int_{z\ge0}\ga e^{\ga z}\,\rmd z P\bigl(G_{\tau(z)-}\in\rmd
r, {\tau(z)} -G_{\tau(z)-}\in\rmd s,X_{\tau(z)}-z\in\rmd\gamma,
\nonumber\hspace*{-15pt}\\
&&\hspace*{155pt} \Xbar
_{\tau(z)-}-X_{\tau(z)-}\in\rmd v, z-\Xbar_{\tau(z)-}\in\rmd y\bigr)\nonumber\hspace*{-15pt}\\
&&\qquad =\cb I(\gamma>0) \int_{z\ge0}\ga e^{\ga z}\,\rmd z I(y\le
z)|V(\rmd r, z-\rmd y)|\whV(\rmd s, \rmd v)\hspace*{-15pt}\\
&&\hspace*{104.4pt}{}\times\Pi_X(\rmd\gamma
+v+y)\nonumber\hspace*{-15pt}\\
&&\qquad\quad{} +\cb\,\rmd_H \int_{z\ge0}\ga
e^{\ga z} \,\rmd z\,\frac{\partial_-}{\partial_- z}V(\rmd t,z)\delta
_0(\rmd s,\rmd\gamma, \rmd v,\rmd y)\nonumber\hspace*{-15pt}\\
&&\qquad =\cb I(\gamma>0)\int_{z\ge y}\ga e^{\ga z}\,\rmd z |V(\rmd r,
z-\rmd y)|\whV(\rmd s, \rmd v)\Pi_X(\rmd\gamma+v+y)\nonumber\hspace*{-15pt}\\
&&\qquad\quad{} +\cb\,\rmd_H \int_{z\ge0}\ga
e^{\ga z}V(\rmd t,\rmd z)\delta_0(\rmd s,\rmd\gamma, \rmd v,\rmd
y).\nonumber\hspace*{-15pt}
\end{eqnarray}
Convolving with the exponential distribution of $\rho$ gives (\ref{102}).

For (\ref{112}), using (\ref{W0}), (\ref{Z}), and independence of
$W$, $Z$ and $\rho$, we have
\begin{eqnarray*}
&&P(G^Z_{\rho-}\in\rmd t, {\rho} -G^Z_{\rho-}\in\rmd s, W_0\in\rmd
\gamma, \Zbar_{\rho-}-Z_{\rho-}\in\rmd\phi, \Zbar_{\rho-}\in
\rmd\theta; W_0> 0)\\
&&\qquad = \cb\ga e^{-\ga\gamma}\,\rmd\gamma P(G^Z_{\rho-}\in\rmd t,
{\rho} -G^Z_{\rho-}\in\rmd s, \Zbar_{\rho-}-Z_{\rho-}\in\rmd
\phi, \Zbar_{\rho-}\in\rmd\theta)\\
&&\qquad = \cb\ga e^{-\ga\gamma}\,\rmd\gamma P\bigl(G^Z_{(t+s)-}\in\rmd t,
Z_{(t+s)-}\in\gt-\rmd\phi, \Zbar_{(t+s)-}\in\rmd\theta\bigr)\\
&&\qquad\quad{}\times \ca
e^{-\ca(t+s)}\,\rmd s\\
&&\qquad = \cb\ga e^{-\ga(\gamma+\phi-\gt)}\,\rmd\gamma e^{\ca(t+s)}\\
&&\qquad\quad{}\times P\bigl(G_{(t+s)-}\in\rmd t, X_{(t+s)-}\in\gt-\rmd\phi, \Xbar
_{(t+s)-}\in\rmd\theta\bigr) \ca e^{-\ca(t+s)}\,\rmd s \\
&&\qquad = \cb\ga e^{-\ga(\gamma+\phi-\gt)}\,\rmd\gamma e^{\ca(t+s)}\\
&&\qquad\quad{}\times P(G_{\rho-}\in\rmd t, \Xbar_{\rho-}\in\rmd\theta, {\rho}
-G_{\rho-}\in\rmd s, \Xbar_{\rho-}-X_{\rho-}\in\rmd\phi)\\
&&\qquad = \cb\ga e^{-\ga(\gamma+\phi-\gt)}\,\rmd\gamma e^{\ca t}
P(G_{\rho-}\in\rmd t, \Xbar_{\rho-}\in\rmd\theta) e^{\ca s}\\
&&\qquad\quad{}\times P({\rho} -G_{\rho-}\in\rmd s, \Xbar_{\rho-}-X_{\rho-}\in\rmd
\phi)
\end{eqnarray*}
by independence of the Wiener--Hopf factors. Further,
\[
e^{\ca t} P(G_{\rho-}\in\rmd t, \Xbar_{\rho-}\in\rmd\theta)=\gk
(\ca, 0)V(\rmd t, \rmd\gt)
\]
and
\[
e^{\ca s} P({\rho} -G_{\rho-}\in\rmd s, \Xbar_{\rho-}-X_{\rho
-}\in\rmd\phi)=\whk(\ca, 0)\whV(\rmd s,\rmd\phi)
\]
as can be seen by taking the Laplace transforms and using (\ref{Vkap})
and (\ref{WHfac}). Equation (\ref{112}) then follows since
$\gk(\ca, 0)\whk(\ca, 0)=\ca$ by (\ref{WH}).
\end{pf}

Theorems \ref{DK10} and \ref{DK11} extend Theorems 10 and 11 in
\cite{DK}.
To see the connection between (\ref{102}) and Theorem 10 of
\cite{DK}, set
\begin{eqnarray*}
m(\rmd t,\rmd y)&=&\int_{r\le t} e^{-\ca r}\,\rmd r \int_{z\ge y} e^{\ga
z}\,\rmd z |V(\rmd t-r,z-\rmd y)|,\\
n(\rmd t,\rmd y)&=&\int_{r\le t} e^{-\ca r}\,\rmd r \int_{z\ge y} e^{\ga
z}\,\rmd z V(\rmd t-r,\rmd z)\delta_0(\rmd y).
\end{eqnarray*}
For any $a>0, b>\ga$,
%
%
\begin{eqnarray}
&&
\int_{t\ge0}\int_{y\ge0}e^{-at-by}m(\rmd t,\rmd y)\nonumber\hspace*{-15pt}\\
&&\qquad=
\int_{t\ge0}\int_{y\ge0}e^{-at-by}\int_{r\le t} e^{-\ca r}\,\rmd r
\int_{z\ge y} e^{\ga z}\,\rmd z |V(\rmd t-r,z-\rmd y)|\nonumber\hspace*{-15pt}\\
&&\qquad=\int_{r\ge0}e^{-\ca r}\,\rmd r\int_{z\ge0}e^{\ga z}\,\rmd z \int
_{t\ge r} \int_{0\le y\le z} e^{-at-by} |V(\rmd t-r,z-\rmd
y)|\nonumber\hspace*{-15pt}\\[-8pt]\\[-8pt]
&&\qquad=\int_{r\ge0}e^{-\ca r}\,\rmd r\int_{z\ge0}e^{\ga z}\,\rmd z \int
_{t\ge0} \int_{0\le y\le z} e^{-a(t+r)-b(z-y)} V(\rmd t,\rmd
y)\nonumber\hspace*{-15pt}\\
&&\qquad=\int_{r\ge0}e^{-(\ca+a) r}\,\rmd r\int_{y\ge0}\int_{t\ge0}
e^{-at+by} V(\rmd t,\rmd y) \int_{z\ge y} e^{-(b-\ga) z}\,\rmd z
\nonumber\\
&&\qquad=\frac1{(\ca+a)\gk(a,-\ga)(b-\ga)}.\nonumber\hspace*{-15pt}
\end{eqnarray}
Similarly,
%
%
\begin{equation}
\int_{t\ge0}\int_{y\ge0}e^{-at-by}n(\rmd t,\rmd y)= \frac1{(\ca
+a)\gk(a,-\ga)}.
\end{equation}
Setting $a=0$ and inverting shows that
%
%
\begin{eqnarray}
\int_{t\ge0}m(\rmd t,\rmd y)&=&\frac{e^{\ga y}\,\rmd y}{\ca\gk(0,-\ga
)}=\frac{e^{\ga y}\,\rmd y}{\cab q},\nonumber\\[-8pt]\\[-8pt]
\int_{t\ge0}n(\rmd t,\rmd y)&=&\frac{\delta_0(\rmd y)}{\ca\gk
(0,-\ga)}=\frac{\delta_0(\rmd y)}{\cab
q}\nonumber
\end{eqnarray}
from (\ref{const}).
Thus after integrating out $t$, (\ref{102}) reduces to
%
%
\begin{eqnarray}\label{1022}
&&P\bigl({\tau_0^W} -G^W_{\too-}\in\rmd s, W_{\tau(0)}\in\rmd\gamma,
\Wbar_{\tau(0)-}-W_{\tau(0)-}\in\rmd v,\nonumber\\
&&\hspace*{153.5pt} -\Wbar_{\tau(0)-}\in\rmd
y; W_0\le0\bigr)\nonumber\\[-8pt]\\[-8pt]
&&\qquad = I(\gamma>0)q^{-1}\ga e^{\ga y}\,\rmd y\whV(\rmd s, \rmd v)\Pi
_X(\rmd\gamma+v+y) \nonumber\\
&&\qquad\quad{} + q^{-1}\ga\,\rmd_H \delta_0(\rmd s, \rmd\gamma
, \rmd v,\rmd y),
\nonumber
\end{eqnarray}
for $\gamma\ge0,s\ge0, v\ge0,y\ge0$. Thus we may conclude that for
$\gamma\ge0,s\ge0,\break v\ge y\ge0$
%
%
\begin{eqnarray}\label{vcon}
&&\lim_{u\to\infty} \Pu\bigl({\tu} -G_{\tu-}\in\rmd s,X_{\tu}-u\in
\rmd\gamma,\nonumber\\
&&\hspace*{59pt} u-X_{\tu-}\in\rmd y, u-\Xbar_{\tu-}\in\rmd
v\bigr)\nonumber\\[-8pt]\\[-8pt]
&&\qquad= I(\gamma>0)q^{-1}\ga e^{\ga y}\,\rmd y\whV(\rmd s, \rmd v-y)\Pi
_X(\rmd\gamma+v)\nonumber\\
&&\qquad\quad{} + q^{-1}\ga\,\rmd_H \delta_0(\rmd s, \rmd\gamma,
\rmd v,\rmd y)
\nonumber
\end{eqnarray}
in the sense of vague convergence.
For $\gamma>0$, this is Doney and Kyprianou's expression in Theorem 10
of \cite{DK}, for the vague limit when $X$ does not creep over the
boundary. The connection between (\ref{112}) and Theorem 11 of
\cite{DK} is similar but easier to see.
It is worth emphasizing that the convergence in Theorem \ref{DK10} is
stronger than the convergence in (\ref{vcon}). In particular,
convergence of the marginals does not follow from the vague convergence
of (\ref{vcon}); indeed, it need not be the case but it does follow
from the weak convergence in Theorem \ref{DK10}. For example, marginal
convergence of the overshoot in (\ref{vcon}) would imply
\begin{eqnarray*}
&&\lim_{u\to\infty} \Pu\bigl(X_{\tu}-u\in\rmd\gamma\bigr)\\
&&\qquad= I(\gamma>0)q^{-1}\ga e^{\ga y}\,\rmd y\int_{y\ge0}\int_{v\ge
y}\int_{s\ge0}\whV(\rmd s, \rmd v-y)\Pi_X(\rmd\gamma+v)\\
&&\qquad\quad{} +
q^{-1}\ga\,\rmd_H \delta_0(\rmd\gamma)\\
&&\qquad= q^{-1}\ga\biggl[\rmd_H\delta_{0}(\rmd\gamma)+ \int_{y\ge0}
e^{\alpha y} \Pi_{ H}(\rmd\gamma+y)\,\rmd y\biggr]
\end{eqnarray*}
by Vigon's equation amicale invers\'ee; see (\ref{DKcor6}) below.
However, by Theorem~\ref{DK10}, in which marginal convergence does
hold, we find that
\begin{eqnarray*}
&&\lim_{u\to\infty} \Pu\bigl(X_{\tu}-u\in\rmd\gamma\bigr)\\
&&\qquad=P\bigl(W_0I(W_0>0)+W_{\tau(0)}I(W_0\le0)\in\rmd\gamma\bigr)\\
&&\qquad=\cb\ga e^{-\alpha\gamma}\,\rmd\gamma
+q^{-1}\ga\biggl[\rmd_H\delta_{0}(\rmd\gamma)+ \int_{y\ge0}
e^{\alpha y} \Pi_{ H}(\rmd\gamma+y)\,\rmd y\biggr]
\end{eqnarray*}
as discussed in Section \ref{s7}.
We will make frequent use of marginal convergence in Theorem \ref
{DK10} in the subsequent sections.


\section{The ruin time}\label{s6}

By taking $f$ constant in the spatial variables in Theorem~\ref{DK10},
we obtain marginal convergence
in the time variables. We begin with the
limiting distribution of the ruin time.
\begin{pf*}{Proof of Theorem \ref{cor1}}
Let $f\dvtx\R\to\R$ be bounded and continuous. Then by (\ref{101a}) [or
(\ref{101b})] with $\gt=0$
\begin{eqnarray*}
\lim_{u\to\infty}E^{(u)}f(\tu)&=&E[f(\rho);W_0> 0]+E\bigl[f\bigl(\rho+\tW
(0)\bigr);W_0\le0\bigr]\\[-1pt]
&=&Ef\bigl(\rho+\tW(0)\bigr),
\end{eqnarray*}
which proves the first equality in (\ref{RT}). Since $\rho+\tW(0)$
has a continuous distribution
%
%
\begin{eqnarray}
P\bigl(\rho+\tW(0)\le t\bigr)
&=& \int_{s\le t}\ca e^{-\ca s}\,\rmd s \cb\int_z \alpha e^{-\alpha
z}\,\rmd z P_z\bigl(\tau(0) < t-s \bigr)\nonumber\\[-1pt]
&=& \int_{s\le t}\ca e^{-\ca s}\,\rmd s \cb\biggl[1+\int_{z>0} \alpha
e^{\alpha z}\,\rmd z P(\Xbar_{t-s}>z) \biggr]\nonumber\\[-8.5pt]\\[-8.5pt]
&=& \cb\int_{s\le t}\ca e^{-\ca s} Ee^{\alpha\Xbar_{t-s}}\,\rmd s
\nonumber\\[-1pt]
&=&\cb E (e^{\alpha\Xbar_{t-\rho}}; \rho\le t)\nonumber
\end{eqnarray}
and the proof is complete.\vspace*{-3pt}
\end{pf*}

Our derivation of the limiting distribution of the ruin time is based
on splitting the distribution at the time of the large jump.
One of the points of distinction between the path decomposition
approach to
studying ruin and that of \cite{DK}, is that in \cite{DK} the split is
at ${G}_{\tau(u)-}$, the time of the last maximum prior to passage over
the boundary. This is a very natural approach given the fluctuation
theory as developed in \cite{bert}, Chapter VI, for example. We now
show how the path decomposition approach can be used to easily derive
the joint limiting distribution of the fluctuation variables
$({G}_{\tau(u)-},\tau(u)-{G}_{\tau(u)-})$ under $P^{(u)}$, thus
extending the results in~\cite{DK}.

Introduce the measures on $[0,\infty)$ given by
%
%
\begin{eqnarray}\label{deldef}\quad
\delta_{\alpha}^V(\rmd t)&=&\int_{\gt\ge0}e^{\alpha\gt}V(\rmd t,
\rmd\gt),\nonumber\\[-1pt]
\delta_{-\alpha}^{\whV}(\rmd s)&=&\int_{\phi\ge0}e^{-\alpha\phi
}\whV(\rmd s,\rmd\phi),
\nonumber\\[-8.5pt]\\[-8.5pt]
K(\rmd s)&=&\int_{z\ge0}(e^{\alpha z} -1)\Pi_{L^{-1},H}(\rmd s,\rmd
z)\nonumber\\[-1pt]
&=&\int_{z\ge0}\ga e^{\ga z}\,\rmd z \pibar_{L^{-1},H}(\rmd
s,z)\nonumber
\end{eqnarray}
and their respective (improper) distribution functions $\delta_{\alpha
}^V(t)$, $\delta_{-\alpha}^{\whV}(s)$ and~$K(s)$, where
\[
\pibar_{L^{-1},H}(\rmd s,z)=\int_{y>z}\pibar_{L^{-1},H}(\rmd s,\rmd y).\vadjust{\goodbreak}
\]

%
\begin{theorem}\label{th1}
Assume (\ref{c1}). Then for all $s,t,\ge0$, we have
%
%
\begin{eqnarray}\label{ttd}
&&\lim_{u \to\infty}
P^{(u)}\bigl({G}_{\tau(u)-}\in\rmd t, \tau(u)-{G}_{\tau(u)-}\in
\rmd s\bigr)\nonumber\\
&&\qquad=\cab\biggl[\delta_\alpha^V (\rmd t)\delta_{-\alpha}^{\whV}(\rmd s)\\
&&\qquad\quad\hspace*{10pt}{}+\bigl(K(\rmd s)+\alpha\,\rmd_H\delta_0(\rmd s)\bigr)
\int_{0\le r\le t}e^{-\ca r}\delta_\alpha^V (\rmd t-r)\,\rmd
r\biggr]\nonumber
\end{eqnarray}
in the sense of weak convergence of probability measures on $[0,\infty)^2$.
\end{theorem}
\begin{pf}
From Theorem \ref{DK10} we have
%
%
\begin{eqnarray}\label{ttd1}
&&\lim_{u \to\infty}
P^{(u)}\bigl({G}_{\tau(u)-}\in\rmd t, \tau(u)-{G}_{\tau(u)-}\in
\rmd s\bigr)\nonumber\\
&&\qquad =P(G^Z_{\rho-}\in\rmd t, {\rho} -G^Z_{\rho
-}\in\rmd s; W_0> 0)\\
&&\qquad\quad{}+ P\bigl(\rho+G^W_{\too-}\in\rmd t, {\tau
_0^W} -G^W_{\too-}\in\rmd s; W_0\le0\bigr).
\nonumber
\end{eqnarray}
Integrating out $\gamma, \theta$ and $\phi$ in (\ref{112})
gives
%
%
\begin{eqnarray}\label{tG1}
&&P(G^Z_{\rho-}\in\rmd t, {\rho} -G^Z_{\rho-}\in\rmd s; W_0>
0)\nonumber\\
&&\qquad= \cab\int_{\gamma>0}\int_{\theta\ge0}\int_{\phi\ge0}\ga
e^{-\ga(\gamma+\phi-\theta)} V(\rmd t,\rmd\theta)\whV(\rmd s,
\rmd\phi)\,\rmd\gamma\\
&&\qquad=\cab\delta_\alpha^V (\rmd t)\delta_{-\alpha}^{\whV}(\rmd
s).\nonumber
\end{eqnarray}
Integrating out $\gamma, y$ and $v$ in the first term of (\ref{102})
gives
%
%
\begin{eqnarray}\label{tG21}
&&P\bigl(\rho+G^W_{\too-}\in\rmd t, {\tau_0^W} -G^W_{\too-}\in\rmd
s,W_{\tau(0)}>0; W_0\le0\bigr)\nonumber\hspace*{-15pt}\\
&&\qquad = \cab\int_{\gamma>0}\int_{y\ge0}\int_{v\ge0} \int
_{r\le t} e^{-\ca r}\,\rmd r \int_{z\ge y}\ga e^{\ga z}\,\rmd z |V(\rmd
t-r,z-\rmd y)|\hspace*{-15pt}\\
&&\qquad\quad\hspace*{155pt}{}\times\whV(\rmd s, \rmd
v)\Pi_X(\rmd\gamma+v+y).\nonumber\hspace*{-15pt}
\end{eqnarray}
By Doney and Kyprianou's extension of Vigon's \'equation amicale
invers\'ee,
it follows that
%
%
\begin{equation}\label{DKcor6}\qquad
\pibar_{L^{-1},H}(\rmd s,y)
=\int_{v\ge0} \whV(\rmd s, \rmd v)\pibar_X^+(v+y),\qquad s\ge0, y\ge0.
\end{equation}
Thus continuing the equalities in (\ref{tG21})
%
%
\begin{eqnarray}\label{tG2}
& = &\cab\int_{y\ge0} \int_{r\le t} e^{-\ca r}\,\rmd r \int
_{z\ge y}\ga e^{\ga z}\,\rmd z |V(\rmd t-r,z-\rmd y)|\pibar
_{L^{-1},H}(\rmd s,y)\nonumber\hspace*{-15pt}\\
& = &\cab\int_{z\ge0} \ga e^{\ga z}\,\rmd z\int_{r\le t} e^{-\ca
r}\,\rmd r \int_{ y\le z} |V(\rmd t-r,z-\rmd y)|\pibar_{L^{-1},H}(\rmd
s,y)\nonumber\hspace*{-15pt}\\
& = &\cab\int_{z\ge0} \ga e^{\ga z}\,\rmd z\int_{r\le t} e^{-\ca
r}\,\rmd r \int_{ y\le z} V(\rmd t-r,\rmd y)\pibar_{L^{-1},H}(\rmd
s,z-y)\hspace*{-15pt}\\
& = &\cab\int_{y\ge0} \int_{r\le t} e^{-\ca r}\,\rmd r V(\rmd
t-r,\rmd y) \int_{z\ge y}\ga e^{\ga z}\,\rmd z \pibar_{L^{-1},H}(\rmd
s,z-y)\nonumber\hspace*{-15pt}\\
& = &\cab\int_{y\ge0} e^{\ga y} \int_{r\le t} e^{-\ca r}\,\rmd r
V(\rmd t-r,\rmd y) \int_{z\ge0}\ga e^{\ga z}\,\rmd z \pibar
_{L^{-1},H}(\rmd s,z)\nonumber\hspace*{-15pt}\\
& = &\cab K(\rmd s)\int_{0\le r\le t}e^{-\ca r}\delta_\alpha^V
(\rmd t-r)\,\rmd
r.\nonumber\hspace*{-15pt}
\end{eqnarray}
Integrating out $\gamma, y$ and $v$ in the second term of (\ref{102})
gives
%
%
\begin{eqnarray}\label{tG3}
&&P\bigl(\rho+G^W_{\too-}\in\rmd t, {\tau_0^W} -G^W_{\too-}\in\rmd s,
W_{\tau(0)}=0; W_0\le0\bigr)\nonumber\\
&&\qquad = \cab\,\rmd_H \int_{r\le t} e^{-\ca r}\,\rmd r\int_{z\ge0}\ga
e^{\ga z} V(\rmd t-r, \rmd z)\delta_0(\rmd s)\\
&&\qquad = \cab\alpha\,\rmd_H\delta_0(\rmd s)\int_0^t e^{-\ca
r}\delta_\alpha^V (\rmd t-r)\,\rmd
r.\nonumber
\end{eqnarray}
Adding the three terms in (\ref{tG1}), (\ref{tG2}) and (\ref{tG3})
gives (\ref{ttd}).
\end{pf}



\section{Overshoots and undershoots}\label{s7}

By taking $f$ constant in the time variables we obtain joint
convergence of overshoots and undershoots.
%
%
\begin{theorem}\label{jtosus}
Assume (\ref{c1}).
Let $f\dvtx\R^2\otimes(\R\cup\{\infty\})\to\R$ be bounded and
jointly continuous. Then for $0\le\gt<\ga$
%
%
\begin{eqnarray}\label{101jo}\qquad
&&\lim_{u\to\infty} \Eu f\bigl(X_{\tu}-u, \Xbar_{\tu-}-X_{\tu-}, \Xbar
_{\tu-}\bigr)e^{-\gt X_{\tu-}I(X_{\tu-}\le0)}\nonumber\\
&&\qquad = E \bigl[f(W_0, \Zbar_{\rho-}-Z_{\rho-}, \Zbar_{\rho-})e^{-\gt
Z_{\rho-}I(Z_{\rho-}\le0)}; W_0> 0\bigr]\\
&&\qquad\quad{} + E \bigl[f\bigl(W_{\tau(0)}, \Wbar_{\tau(0)-}-W_{\tau(0)-}, \infty\bigr);
W_0\le
0\bigr]\nonumber
\end{eqnarray}
and
%
%
\begin{eqnarray}\label{101ju}
&&\lim_{u\to\infty} \Eu f\bigl(X_{\tu}-u, \Xbar_{\tu-}-X_{\tu-},
u-\Xbar_{\tu-}\bigr)e^{-\gt X_{\tu-}I(X_{\tu-}\le0)}\hspace*{-35pt}\nonumber\\
&&\qquad = E \bigl[f( W_0, \Zbar_{\rho-}-Z_{\rho-}, \infty)e^{-\gt Z_{\rho
-}I(Z_{\rho-}\le0)}; W_0> 0\bigr]\hspace*{-35pt}\\
&&\qquad\quad{} + E \bigl[f\bigl( W_{\tau(0)}, \Wbar_{\tau(0)-}-W_{\tau(0)-},- \Wbar
_{\tau(0)-}\bigr); W_0\le
0\bigr].\nonumber\hspace*{-35pt}
\end{eqnarray}
For $\gamma> 0$, $\theta\ge0, \phi\ge0$
%
%
\begin{eqnarray}\label{112jo}
&&
P( W_0\in\rmd\gamma, \Zbar_{\rho-}-Z_{\rho-}\in\rmd\phi, \Zbar
_{\rho-}\in\rmd\theta; W_0> 0)\nonumber\\[-8pt]\\[-8pt]
&&\qquad=
\cab\ga e^{-\ga(\gamma+\phi-\gt)} V(\rmd\theta)\whV(\rmd\phi
)\,\rmd\gamma.\nonumber
\end{eqnarray}
For $\gamma\ge0, v\ge0, y\ge0$
%
%
\begin{eqnarray}\label{102jo}
&&P\bigl( W_{\tau(0)}\in\rmd\gamma, \Wbar_{\tau(0)-}-W_{\tau(0)-}\in
\rmd v, -\Wbar_{\tau(0)-}\in\rmd y; W_0\le0\bigr)\nonumber\\
&&\qquad = I(\gamma>0)q^{-1}\ga e^{\ga y}\,\rmd y\whV( \rmd v)\Pi
_X(\rmd\gamma+v+y) \\
&&\qquad\quad{ }+ q^{-1}\ga\,\rmd_H \delta_0(\rmd\gamma, \rmd
v,\rmd
y).\nonumber
\end{eqnarray}
\end{theorem}
\begin{pf}
The result follows immediately from (\ref{101a}), (\ref{101b}),
(\ref{112})\break and~(\ref{1022}).
\end{pf}

Theorem \ref{jtosus}
contains all results we know of
in the literature on convergence of individual overshoots and
undershoots, under a convolution equivalent assumption.
The only marginal limiting distribution in Theorem
\ref{jtosus} which is proper is that of the overshoot, and this is
given by $W_0I(W_0>0)+W_{\tau(0)}I(W_0\le0)$.
An easy calculation from (\ref{112jo}) and (\ref{102jo}), using
(\ref{DKcor6}), gives the following.

\textit{Overshoot}. Assume (\ref{c1}). Then for
$\gamma\ge0$
%
%
\begin{eqnarray}\label{o}
&&\Pu\bigl(X_{\tu}-u\in\rmd\gamma\bigr)\nonumber\\
&&\qquad \to\cb\ga e^{-\alpha\gamma}\,\rmd
\gamma\\
&&\qquad\quad{}+q^{-1}\ga\biggl[\rmd_H\delta_{0}(\rmd\gamma)+ \int_{y\ge0}
e^{\alpha y} \Pi_{ H}(\rmd\gamma+y)\,\rmd y\biggr].\nonumber
\end{eqnarray}
Observe that the limiting distribution has mass $\alpha\,\rmd
_{H}q^{-1}$ at the origin, and for $x>0$
%
%
\begin{equation}\label{pm}\qquad
\Pu\bigl(X_{\tu}-u>x\bigr)\to\cb e^{-\alpha x}+q^{-1}e^{-\alpha x}
\int_{y>x}(e^{\alpha y}- e^{\alpha x})\Pi_{ H}(\rmd y).
\end{equation}

This is the form of the limiting distribution given in \cite{kkm} and
\cite{DK}.
In~\cite{DK}, it is indicated that the limiting distribution on
$(0,\infty)$ arises as a consequence of either an arbitrarily large
jump from a finite
position after a~finite time, or a~finite jump from a finite distance
relative to the boundary after an arbitrarily
large time. This is not quite correct. From the path decomposition, the
latter component of the limiting distribution arises
as a~consequence of a~large jump from a finite position to within a
finite distance of the boundary after a finite time,
followed by a finite jump a finite time later. The atom at $0$ in the
limiting distribution is a consequence of creeping across
the boundary when the large jump undershoots the boundary.\looseness=1

The other marginal limits in Theorem
\ref{jtosus} are improper, thus in each instance below, convergence is
in the vague sense with remaining mass escaping to $+\infty$.
We leave the calculations to the reader.

\textit{Undershoots}. Assume (\ref{c1}). Then for
$x\ge0$
%
%
\begin{eqnarray}\label{u}
&&P^{(u)}\bigl(u-X_ {{\tu}-}\in\rmd x\bigr)\nonumber\\[-8pt]\\[-8pt]
&&\qquad\to
q^{-1} \alpha\,\rmd_H\delta_{0}(\rmd x)
+q^{-1} \alpha e^{\alpha x} \pibar_X(x)\,\rmd x\int_{0\le v\le
x}e^{-\ga v}\whV(\rmd v),\nonumber
\end{eqnarray}
%
while for $y\ge0$
%
%
\begin{equation}\label{mu}\quad
P^{(u)}\bigl(u-\Xbar_ {{\tu}-}\in\rmd y\bigr)\to
q^{-1} \alpha\,\rmd_H\delta_{0}(\rmd y)
+q^{-1} \alpha
e^{\alpha y}\pibar_{H}(y)\,\rmd y.
\end{equation}

%
\begin{remark}
An alternative formulation of (\ref{u})
appears in Theorem~3.2 of~\cite{pama}.
Statement (\ref{mu}) corrects an
oversight in Theorem 3.3 of \cite{pama}. The first term in
(\ref{mu}), representing possible
mass at 0 if creeping over the boundary occurs, was omitted.
\end{remark}
\textit{Positions prior to overshoot}. Assume (\ref{c1}). Then for
$\gz\in(-\infty, \infty)$
%
%
\begin{equation}
\label{up}
P^{(u)}\bigl(X_ {{\tu}-}\in\rmd\gz\bigr)\to\cab e^{\alpha\gz}V_X(\rmd\gz),
\end{equation}
%
where $V_X$ is the potential measure of $X$, while for $\gt\ge0$
\[
P^{(u)}\bigl(\Xbar_ {{\tu}-}\in\rmd\gt\bigr)\to q\cb^2 e^{\alpha\gt
}V(\rmd\gt).
\]


\section{Laplace transforms and penalty functions}\label{s8}

Expected discounted penal\-ty functions (EDPFs) were introduced into risk
theory by Gerber and Shiu~\cite{GS}.
As an example, consider
%
%
\begin{equation}\label{ltLT}
\Eu e^{-\nu G_{\tu-}-\gz(\tau_u-G_{\tu-})-\eta(X_{\tau_u}-u)
-\lambda(u-\Xbar_{\tau_u-})},
\end{equation}
where $\nu\ge0$, $\gz\ge0$, $\eta>-\ga$, $\lambda\ge0$. In this
case, penalization is more severe when the shortfall at ruin is greater
(if $\eta<0$), but this is moderated by a~later occurrence of ruin or
by a larger minimum surplus prior to ruin.
Among other things, EDPFs provide a natural approach to studying
solvency requirements, and more generally, to valuing cash-flows
related to first passage over a barrier; see, for example, the
discussion in Biffis and Morales \cite{BM}.
In this section we use our previous results to calculate the limit, as
$u\to\infty$, of~(\ref{ltLT}) and other related EDPFs and Laplace transforms.

If $\eta\ge0$, then the limit in (\ref{ltLT}) can be found by using
Theorem \ref{DK10}.
To include the case $-\ga<\eta<0$, it will suffice, by
uniform integrability, to show that
%
%
\begin{equation}\label{mgf1}
\limsup_{u\to\infty} \Eu e^{-\eta(X_{\tu}-u)}<\infty,\qquad \eta
>-\ga.
\end{equation}
A stronger version of (\ref{mgf1}) is in Park and Maller \cite{pama}.
Since our weaker version is easy to prove, we give a direct proof that
does not involve delicate estimation of convolution equivalent
integrals as in \cite{pama}. Combined with convergence of the
overshoot, this weaker result is in fact equivalent to Park and
Maller's {a priori} stronger result on convergence of the mgf of the overshoot.

%
\begin{lemma}\label{EG}
Let $F$ and $G$ be distribution functions with $F(0-)=\break G(0-)=0$, $F\in
\mathcal{S}^{(\alpha)}$ and
%
%
\begin{equation}\label{GF}
\limsup_{u\to\infty}\frac{\overline{G}(u)}{\overline
{F}(u)}<\infty.
\end{equation}
Then
%
%
\begin{equation}\label{EG1}
\limsup_{u\to\infty} \int\frac{\overline{F}(u-y)}{\overline
{F}(u)}G(\rmd y)<\infty.
\end{equation}
\end{lemma}
\begin{pf}
Statement (\ref{GF}) implies $\sup_{u}\overline{G}(u)/\overline{F}(u)
\le C$ for some $C<\infty$,
so the lemma follows easily from (\ref{S2}) since
\begin{eqnarray*}
\int{\overline{F}(u-y)}G(\rmd y)
&=& \int{\overline{G}(u-y)}F(\rmd y)\\
&\le& C\int{\overline{F}(u-y)}F(\rmd y)\\
&=& C{\overline{F}{}^{*2}(u)}.
\end{eqnarray*}
\upqed\end{pf}

In the following lemma, $C$ denotes an unimportant constant which may
change in value from one usage to the next.
%
%
\begin{lemma}\label{obdd}
For any $\eta>-\ga$,
%
%
\begin{equation}\label{mgf}
\limsup_{u\to\infty} \Eu e^{-\eta(X_{\tu}-u)}<\infty.
\end{equation}
\end{lemma}

\begin{pf}
Let
$
T(u)=\inf\{t\dvtx H_t>u\}.
$
Then $\tau(u)=L^{-1}_{T(u)}$ and $X_{\tau(u)}=H_{T(u)}$. Hence,
applying the killed version of Proposition III.2 of \cite{bert} (see
\cite{kypbook}, Theorem~5.6), for $x\ge0$
%
%
\begin{eqnarray}\label{kos}\qquad
\Pu\bigl(X_{\tu}-u>x\bigr)&=&\frac{P(H_{T(u)}-u>x,T(u)<\infty)}{P(\tu<\infty
)}\nonumber\\[-8pt]\\[-8pt]
&=&\frac{\pibar_H(u)}{P(\tu<\infty)}\int_{0\le y\le u} \frac
{\pibar_H(u-y+x)}{\pibar_H(u)}V(\rmd
y).\nonumber
\end{eqnarray}
Fix $\veps>0$ so that $\ga-\veps+\eta>0$. Applying (\ref{Pot}),
%
%
\begin{eqnarray}\label{obdd1}
&&\int_{0\le y\le u-1} \frac{\pibar_H(u-y+x)}{\pibar_H(u)}V(\rmd y)\nonumber\\
&&\qquad\le Ae^{-(\ga-\veps)x}\int_{0\le y\le u-1} \frac{\pibar
_H(u-y)}{\pibar_H(u)}V(\rmd y)\\
&&\qquad\le ACe^{-(\ga-\veps)x},\nonumber
\end{eqnarray}
if $u\ge2$, since by (\ref{clim}), (\ref{XHlim}) and (\ref{HSga}),
Lemma \ref{EG} may be applied to the distributions
$F(\rmd y)=I(y> 1)\Pi_H(\rmd y)/\pibar_H(1)$ and $G(\rmd y)=V(\rmd
y)/V(\infty)$. On the other hand,
%
%
\begin{equation}\label{obdd2}
\int_{u-1\le y\le u} \frac{\pibar_H(u-y+x)}{\pibar_H(u)}V(\rmd
y)\le\pibar_H(x)\frac{\Vbar(u-1)}{\pibar_H(u)}\le Ce^{-\ga x}
\end{equation}
as $\Pi_H\in\mathcal{S}^{(\alpha)}$.
Since the ratio in front of the integral in (\ref{kos}) is bounded
by~(\ref{clim})--(\ref{XHlim}), the result follows from (\ref{obdd1})
and (\ref{obdd2}).
\end{pf}

As preparation for calculating the limit of (\ref{ltLT}) we need the
following proposition.
%
%
\begin{prop}\label{LTall}
Let $\nu\ge0$, $\gz\ge0$, $\eta>-\ga$, $\lambda\ge0$.
Then
%
%
\begin{eqnarray}\label{LTZ}
&&
E \bigl[e^{-\nu G^Z_{\rho-}-\gz(\rho-G^Z_{\rho-})-\eta W_0-\lambda\Zbar
_{\rho-}};W_0> 0\bigr]\nonumber\\[-8pt]\\[-8pt]
&&\qquad=\frac{\cab\ga\gk(\gz,-\ga)}
{(\ga+\eta)(\gz+\ca)\kappa(\nu,\lambda-\ga)}.\nonumber
\end{eqnarray}
If, in addition, $\lambda\ne\ga+\eta$, then
%
%
\begin{eqnarray}\label{LTW}
&& E \bigl[e^{-\nu(\rho+G^W_{\too-})-\gz(\tau_0^W-G^W_{\too-})-\eta
W_{\too}+\lambda\Wbar_{\too-}};W_0\le0\bigr]\nonumber\\[-8pt]\\[-8pt]
&&\qquad =\frac{\cab\ga[\kappa
(\gz,\lambda-\ga)-\kappa(\gz,\eta)]}
{(\ca+\nu)(\lambda-\ga-\eta)\kappa(\nu,-\ga)}.\nonumber
\end{eqnarray}
\end{prop}
\begin{pf}
Fix $\nu\ge0$, $\gz\ge0$, $\eta>-\ga$, $\lambda\ge0$. Then by
(\ref{112})
\begin{eqnarray*}
&&E \bigl[e^{-\nu G^Z_{\rho-}-\gz(\rho-G^Z_{\rho-})-\eta W_0-\lambda
\Zbar_{\rho-}};W_0> 0\bigr]\\
&&\qquad = \cab\ga\int_{t\ge0}\int_{s\ge0} \int_{ \gamma>0}\int
_{\gt\ge0}\int_{\phi\ge0}
e^{-\nu t-\gz s-\eta\gamma-\lambda\gt}e^{-\ga(\gamma+\phi-\gt)}\\
&&\hspace*{155pt}{}\times V(\rmd t,\rmd\theta)\whV(\rmd s, \rmd\phi)\,\rmd\gamma\\
&&\qquad = \cab\ga\int_{ \gamma>0}e^{-(\ga+\eta) \gamma} \,\rmd
\gamma\int_{t\ge0}\int_{\gt\ge0}
e^{-\nu t-(\lambda-\ga) \gt} V(\rmd t,\rmd\theta) \\
&&\qquad\quad{}\times\int_{s\ge0}
\int_{\phi\ge0}e^{-\gz s-\ga\phi}\whV(\rmd s,\rmd\phi)\\
&&\qquad =\frac{\cab\ga}
{(\ga+\eta)\kappa(\nu,\lambda-\ga)\whk(\gz,\ga)},
\end{eqnarray*}
since $\kappa(\nu,\lambda-\ga)>0$ by (\ref{kpos}) and $\whk(\gz
,\ga)>0$ trivially.
Thus (\ref{LTZ}) follows from (\ref{WH1}) and (\ref{const}).

Now assume $\lambda\ne\ga+\eta$, then by (\ref{defW})
\begin{eqnarray*}
&&E \bigl[e^{-\nu(\rho+G^W_{\too-})-\gz(\tau_0^W-G^W_{\too-})-\eta
W_{\too}+\lambda\Wbar_{\too-}};W_0\le0\bigr]\\[-3pt]
&&\qquad = \cb\ga\int_{z\le0} e^{-\ga z} \,\rmd z E_z \bigl[e^{-\nu(\rho
+G_{\too-})-\gz(\too-G_{\too-})-\eta X_{\too}+\lambda\Xbar_{\too
-}};\\[-3pt]
&&\qquad\hspace*{248pt}\too<\infty\bigr]\\[-3pt]
&&\qquad = \frac{\cab\ga}{\ca+\nu} \int_{z> 0} e^{\ga z} \,\rmd z E
\bigl[e^{-\nu G_{\tau(z)-}-\gz(\tau(z)-G_{\tau(z)-})-\eta(X_{\tau
(z)}-z)-\lambda(z-\Xbar_{\tau(z)-})};\\[-3pt]
&&\hspace*{288pt}\tau(z)<\infty\bigr]\\[-3pt]
&&\qquad =\frac{\cab\ga[\kappa(\gz,\lambda-\ga)-\kappa(\gz
,\eta)]}
{(\ca+\nu)(\lambda-\ga-\eta)\kappa(\nu,-\ga)}
\end{eqnarray*}
by the extension of the second factorization identity in Theorem 3.5 of
\cite{GM1}.
\end{pf}

We are now ready to calculate the limit of (\ref{ltLT}) and a related
penalty function.
%
%
\begin{theorem}\label{LTlt}
Fix $\nu\ge0$, $\gz\ge0$, $\eta> -\ga$, $\lambda> 0$.
Then
%
%
\begin{eqnarray}\label{LTlt0}
&&\lim_{u\to\infty} \Eu e^{-\nu G_{\tu-}-\gz(\tau_u-G_{\tu-})-\eta
(X_{\tau(u)}-u)-\lambda\Xbar_{\tau(u)-}}\nonumber\\[-8pt]\\[-8pt]
&&\qquad=\frac{\cab\ga\gk(\gz,-\ga)}
{(\ga+\eta)(\gz+\ca)\kappa(\nu,\lambda-\ga)}.\nonumber
\end{eqnarray}
If, in addition, $\lambda\ne\ga+\eta$, then
%
%
\begin{eqnarray}\label{LTltu}
&&\lim_{u\to\infty} \Eu e^{-\nu G_{\tu-}-\gz(\tau(u)-G_{\tu
-})-\eta(X_{\tau(u)}-u)-\lambda(u-\Xbar_{\tau(u)-})}\nonumber\\[-8pt]\\[-8pt]
&&\qquad=\frac{\cab\ga[\kappa(\gz,\lambda-\ga)-\kappa(\gz,\eta)]}
{(\ca+\nu)(\lambda-\ga-\eta)\kappa(\nu,-\ga)}.\nonumber
\end{eqnarray}
\end{theorem}
\begin{pf}
Since (\ref{LTlt0}) and (\ref{LTltu}) follow in a similar manner from
(\ref{LTZ}) and (\ref{LTW}), respectively, we only prove (\ref{LTlt0}).
Let
\[
g(t,s,\gamma, y)=e^{-\nu t-\gz s -\eta\gamma-\lambda y}.
\]
By (\ref{101aa}), (\ref{mgf}) and uniform integrability
\begin{eqnarray*}
&&\lim_{u\to\infty} \Eu g\bigl(G_{\tu-}, {\tu} -G_{\tu-},X_{\tu}-u,
\Xbar_{\tu-}\bigr)\\
&&\qquad = E [g(G^Z_{\rho-},{\rho} -G^Z_{\rho-}, W_0, \Zbar_{\rho
-}); W_0> 0]\\
&&\qquad\quad{} + E \bigl[g\bigl(\rho+G^W_{\too-},{\tau_0^W} -G^W_{\too-},
W_{\tau(0)}, \infty\bigr); W_0\le0\bigr]\\
&&\qquad = E \bigl[e^{-\nu G^Z_{\rho-}-\gz({\rho} -G^Z_{\rho-})-\eta
W_0+\lambda\Zbar_{\rho-}};W_0> 0\bigr],
\end{eqnarray*}
since $\lambda>0$.
Thus (\ref{LTlt0}) follows from (\ref{LTZ}).\vadjust{\goodbreak}
\end{pf}

Setting $\eta=\nu=\gz=0$ in (\ref{LTlt0}) gives
%
%
\begin{eqnarray}\label{PM1}
\lim_{u\to\infty}e^{-\lambda u} \Eu e^{\lambda(u-\Xbar_{\tau
(u)-})}&=&\lim_{u\to\infty} \Eu
e^{-\lambda\Xbar_{\tau(u)-}}\nonumber\\[-8.5pt]\\[-8.5pt]
&=&\frac{\cb\gk(0,-\ga)}{\kappa(0,\lambda-\ga)}.\nonumber
\end{eqnarray}
This gives a transparent explanation of the mgf result in Theorem 3.3
of~\cite{pama}, and extends it to all $\gl> 0$. Note that letting
$\gl\downarrow0$ in the final expression of~(\ref{PM1}), reflects
that in the limit, $\Xbar_{\tau(u)-}$ has mass $1-\cb$ at infinity
under~$\Pu$. Similarly setting $\eta=\nu=\gz=0$ in (\ref{LTltu})
gives the growth in the mgf of $\Xbar_{\tau(u)-}$ as measured from
the origin; for every $\gl> 0$
%
%
\begin{eqnarray}\label{PM11}
\lim_{u\to\infty}e^{-\lambda u} \Eu e^{\lambda\Xbar_{\tau(u)-}}
&=&\frac{\cb\ga[\gk(0,\gl-\ga)-q]}{(\lambda-\ga)\kappa(0,-\ga
)}\nonumber\\[-8.5pt]\\[-8.5pt]
&=&\frac{\ga[\gk(0,\gl-\ga)-q]}{(\lambda-\ga)q}.\nonumber
\end{eqnarray}
In this case, letting $\gl\downarrow0$ reflects that in the limit,
$u-\Xbar_{\tau(u)-}$ has mass $\cb$ at infinity under $\Pu$.

Observe that (\ref{LTlt0}) and (\ref{LTltu}) are both false when $\gl
=0$, as can be seen from (\ref{PM1}) and (\ref{PM11}). In this case,
the limit is obtained by adding the corresponding expressions in (\ref
{LTZ}) and (\ref{LTW}).
%
%
\begin{theorem}\label{ltjt}
Fix $\nu\ge0$, $\gz\ge0$, $\eta> -\ga$,
then
%
%
\begin{eqnarray}\label{ltjt1}
&&\lim_{u\to\infty} \Eu e^{-\nu G_{\tu-}-\gz(\tau(u)-G_{\tu
-})-\eta(X_{\tau(u)}-u)}\nonumber\\[-8.5pt]\\[-8.5pt]
&&\qquad =\frac{\cab\ga}
{(\ga+\eta)\kappa(\nu,-\ga)}\biggl[\frac{\gk(\gz,-\ga)}{\ca
+\gz}+\frac
{\gk(\gz,\eta)-\gk(\gz,-\ga)}{\ca+\nu}\biggr].\nonumber
\end{eqnarray}
\end{theorem}
\begin{pf}
By Theorem \ref{DK10} and (\ref{mgf}),
\begin{eqnarray*}
&&\lim_{u\to\infty} \Eu e^{-\nu G_{\tu-}-\gz(\tau(u)-G_{\tu
-})-\eta(X_{\tau(u)}-u)}\\[-1pt]
&&\qquad=E \bigl[e^{-\nu G^Z_{\rho-}-\gz(\rho-G^Z_{\rho-})-\eta W_0};W_0> 0\bigr]\\[-1pt]
&&\quad\qquad{}+
E \bigl[e^{-\nu(\rho+G^W_{\too-})-\gz(\tau_0^W-G^W_{\too-})-\eta
W_{\too}};W_0\le0\bigr].
\end{eqnarray*}
The result now follows by setting $\gl=0$ in (\ref{LTZ}) and (\ref
{LTW}) and adding.
\end{pf}

As a special case of (\ref{ltjt1}), with $\nu=\gz$, we obtain the
limit of the joint transform of the overshoot and ruin time,
%
%
\begin{equation}\label{LTto}
\lim_{u\to\infty} \Eu e^{-\gz\tau(u)-\eta(X_{\tau(u)}-u)}
=\frac{\cab\ga\gk(\gz,\eta)}
{(\ga+\eta)(\ca+\gz)\kappa(\gz,-\ga)}.
\end{equation}
Setting $\gz=0$ in (\ref{LTto})
evaluates the limit in (\ref{mgf1}). With $\eta=0$, (\ref{LTto})
reflects the description of the limiting distribution in Theorem
\ref{cor1}.\vadjust{\goodbreak}

We now briefly describe an application of the EDPF in (\ref{LTto})
when $\eta>0$. Fix $\gz\ge0$ and choose $\eta=\eta(\gz)$
so that $e^{-\gz t-\eta X_{t}}$ is a martingale.
In actuarial terms, $\eta$ is a solution to
\textit{Lundberg's fundamental equation}
(see, e.g., Gerber-Shiu \cite{GS}, page~51).
To see that such an $\eta$ exists and is unique in our setup,
first observe that by (\ref{WH}), this is equivalent to
%
%
\begin{equation}\label{LFE}
\gk(\gz,\eta)\whk(\gz,-\eta)=0.
\end{equation}
Now for $\gz\ge0$,
\[
e^{-\gk(\gz,\eta)}=e^{-q}E e^{-\gz\mathcal{L}^{-1}_1 -\eta
\mathcal{H}_1}\le e^{-q}Ee^{-\eta\mathcal{H}_1}=
\cases{
<1, &\quad $\eta\ge-\ga$,\cr
=\infty, &\quad $\eta< -\ga$,}
\]
by Proposition 5.1 of \cite{kkm}.
Thus in order that (\ref{LFE}) holds, it must be that
$\whk(\gz,-\eta)=0$. Since $\whk(\gz,0)\ge0$ and
$\whk(\gz,-\eta)\downarrow-\infty$ as $\eta\uparrow\infty$,
this equation has a unique solution $\eta\ge0$. Then by (\ref
{LTto}), if $\gz>0$ and $\eta=\eta(\gz)$,
\[
\Eu\frac{e^{-\gz\tu}(1 -e^{-\eta(X_{\tau(u)}-u)})}{\gz}
\to\frac{\cab\ga}{\gz(\ca+\gz)\gk(\gz,-\ga)}\biggl(\frac{\gk
(\gz,0)}{\ga}-\frac{\gk(\gz,\eta)}{\ga+\eta}\biggr).
\]
In the spectrally positive case, Gerber and Shiu \cite{GS} interpret
this in terms of the expected present value of a deferred continuous
annuity at a rate of 1 per unit time,
starting at the time of ruin and ending as soon as the shortfall
returns to zero.

The standard form of the EDPFs introduced by Gerber and Shiu is
%
%
\begin{equation}\label{GSpf}
\Eu\bigl[e^{-\gz\tu}g\bigl( X_{\tu}-u, u-X_{\tu-}\bigr)\bigr]
\end{equation}
for suitably chosen functions $g$.
We have chosen to formulate the results in this section in terms of exponential
penalty functions using the undershoot of the maximum $u-\Xbar_{\tu
-}$ instead of $u-X_{\tu-}$. It is clear that more general penalty
functions could have been used, and the resulting limits could then be
found using Theorems \ref{DK10} and \ref{DK11}.
For the Gerber and Shiu penalty function in (\ref{GSpf}), under the
appropriate conditions on $g$ so that Theorem \ref{DK10} applies, we have
%
%
\begin{eqnarray}\label{GSpf1}
&&\lim_{u\to\infty}\Eu\bigl[e^{-\gz\tu}g\bigl( X_{\tu}-u, u-X_{\tu-}\bigr)\bigr] \nonumber\\
&&\qquad = E [e^{-\gz\rho}g(W_0, \infty); W_0> 0] \\
&&\qquad\quad{}+ E \bigl[e^{-\gz(\rho+{\tau
_0^W})}g\bigl(W_{\tau(0)},-W_{\tau(0)-}\bigr); W_0\le
0\bigr].\nonumber
\end{eqnarray}
A natural example would be
%
%
\begin{eqnarray}\label{GSeg}
&&
\lim_{u\to\infty}\Eu e^{-\gz\tu-\eta(X_{\tu}-u)-\lambda
(u-X_{\tu-})}\nonumber\\[-8pt]\\[-8pt]
&&\qquad=Ee^{-\gz{\rho}}E\bigl[e^{-\gz{\tau_0^W}-\eta W_{\tau(0)}+\lambda
W_{\tau(0)-}};W_0\le0\bigr]\nonumber
\end{eqnarray}
for $\gz\ge0$, $\eta>-\ga$ and $\lambda> 0$. The limit can then be
calculated using Theorem~\ref{DK11}, although the resulting expression
obtained is not as simple as those obtained in Theorem \ref{LTlt}.
Quite different behavior occurs if $\gl<0$ in (\ref{GSeg}).
%
%
\begin{prop}\label{GSeg0}
Let $\nu\ge0$, $\gz\ge0$, $\eta>-\ga$, $0<\gt<\ga$ and assume
that $\gt-\eta<\ga$. Then
%
%
\begin{eqnarray}\label{GSeg1}
&&\lim_{u\to\infty}\Eu e^{-\nu G_{\tu-}-\gz(\tau(u)-G_{\tu-})-\eta
(X_{\tu}-u)-\gt X_{\tu-}}\nonumber\\[-8pt]\\[-8pt]
&&\qquad=\frac{\cab\ga}{(\alpha+\eta)\gk(\nu,
\gt-\ga)\whk(\gz, \ga-\gt)}.\nonumber
\end{eqnarray}
\end{prop}
\begin{pf}
With $\eta$ and $\gt$ as above, we first observe that
%
%
\begin{equation}\label{ouui}
\limsup_{u\to\infty}\Eu e^{-\eta(X_{\tu}-u)-\gt X_{\tu-}I(X_{\tu
-}\le0)}<\infty.
\end{equation}
This follows immediately from (\ref{undbd}) if $\eta\ge0$;
consequently we may assume $-\ga<\eta<0$. By separately considering
the cases
$X_{\tu}-u> |X_{\tu-}|$ and $X_{\tu}-u\le|X_{\tu-}|$, one finds
\[
e^{-\eta(X_{\tu}-u)-\gt X_{\tu-}I(X_{\tu-}\le0)}\le e^{(\gt-\eta
)(X_{\tu}-u)}+e^{-(\gt-\eta)X_{\tu-}I(X_{\tu-}\le0)},
\]
and so (\ref{ouui}) again follows from (\ref{undbd}) and (\ref
{mgf}), since $\gt-\eta<\ga$.
Hence, $e^{-\eta(X_{\tu}-u)-\gt X_{\tu-}}$ is uniformly integrable
if $\eta>-\ga$, $0<\gt<\ga$ and $\gt-\eta<\ga$.
Thus by Theorems \ref{DK10} and \ref{DK11}
\begin{eqnarray*}
&&\lim_{u\to\infty}\Eu e^{-\nu G_{\tu-}-\gz(\tau(u)-G_{\tu
-})-\eta(X_{\tu}-u)-\gt X_{\tu-}}\\
&&\qquad =E\bigl[e^{-\nu G^Z_{\rho-}-\gz(\rho-G^Z_{\rho-})-\eta W_0 -\gt
Z_{\rho-} };W_0>0\bigr]\\
&&\qquad =\cab\ga\int_{t\ge0}\int_{s\ge0}\int_{\gamma>0}\int
_{\phi\ge0}\int_{\xi\ge-\phi}
e^{-\nu t-\gz s-\eta\gamma-\gt\xi}\\
&&\hspace*{165pt}{}\times e^{-\ga(\gamma-\xi)} V(\rmd t,\phi+ \rmd\xi)\whV(\rmd s,\rmd
\phi)\,\rmd\gamma\\
&&\qquad =\frac{\cab\ga}{\ga+\eta}\int_{t\ge0}\int_{\xi\ge0}
e^{-\nu t+(\ga-\gt)\xi}V(\rmd t,\rmd\xi)\\
&&\qquad\hspace*{57pt}{}\times\int_{s\ge0}\int_{\phi\ge0} e^{-\gz s-(\ga-\gt)\phi}V(\rmd
t,\rmd\theta)\whV(\rmd s, \rmd\phi),
\end{eqnarray*}
which gives (\ref{GSeg1}).
\end{pf}

Setting $\nu=\gz$, using (\ref{WH1}) and rewriting (\ref{GSeg1}) in
terms of the undershoot gives
%
%
\begin{eqnarray}\label{GSeg2}
&&\lim_{u\to\infty}e^{-\gt u}\Eu e^{-\gz\tu-\eta(X_{\tu}-u)+\gt
(u-X_{\tu-})}\nonumber\\[-8pt]\\[-8pt]
&&\qquad= \frac{\cab\ga}{(\alpha+\eta)(\gz-\Psi(i(\gt-\ga)))}.\nonumber
\end{eqnarray}
The special case of (\ref{GSeg2}) with $\gz=\eta=0$ is given in
Theorem 3.2 of \cite{pama}.
Results related to (\ref{GSeg2}) for the case of a Cram\'er--Lundberg
model with bounded claims density can be found in Corollary 3.2 of Tang
and Wei \cite{TW}.
When $\gt=0$, (\ref{GSeg1}) fails just as (\ref{LTlt0}) fails when
$\gl=0$. Observe though that letting $\gt\downarrow0$ on the RHS of
(\ref{GSeg1})
and $\gl\downarrow0$ on the RHS of (\ref{LTlt0}) results in the same
limit, as one would expect.







\printaddresses

\end{document}